\def\ifundefined#1#2{\expandafter\ifx\csname#1\endcsname\relax\input #2\fi}
\def\Ma{{\cal\char'101}}
\def\Mb{{\cal\char'102}}
\def\Mc{{\cal\char'103}}
\def\Md{{\cal\char'104}}

\def\Mg{{\cal\char'107}}

\def\Mi{{\cal\char'111}}
\def\Mj{{\cal\char'112}}

\def\Mm{{\cal\char'115}}
\def\Mn{{\cal\char'116}}
\def\Mo{{\cal\char'117}}

\def\Mr{{\cal\char'122}}
\def\Ms{{\cal\char'123}}
\def\Mt{{\cal\char'124}}
\def\Mu{{\cal\char'125}}
\def\Mv{{\cal\char'126}}

\def\Mx{{\cal\char'130}}

\def\pmb#1{\setbox0=\hbox{$#1$}       
     \kern-.025em\copy0\kern-\wd0
     \kern.05em\copy0\kern-\wd0
     \kern-.025em\box0}


\def\endproof{$\hfill \square$}




\def\Cross{\bigm| \kern-5.5pt \not \ \, }
\def\cross{\mid \kern-5.0pt \not \ \, }             
\def\notto{\hbox{$~\rightarrow~\kern-1.5em\hbox{/}\ \ $}}

\def\al{\alpha}
\def\be{\beta}

\def\vph{\varphi}

\def\sg{\sigma}

\hyphenation{math-ema-ticians}
\hyphenation{pa-ra-meters}
\hyphenation{pa-ra-meter}
\hyphenation{lem-ma}
\hyphenation{lem-mas}
\hyphenation{to-po-logy}
\hyphenation{to-po-logies}
\hyphenation{homo-logy}
\hyphenation{homo-mor-phy}

\def\nSigma{\Sigma \kern-8.3pt \bigm|\,}

\def\got#1{\hbox{\teneuler #1}}

\font\teneufm=eufm10
\font\eighteufm=eufm8
\font\fiveeufm=eufm5

\newfam\eufam
\textfont\eufam=\teneufm
\scriptfont\eufam=\eighteufm
\scriptscriptfont\eufam=\fiveeufm

\def\got{\fam=\eufam\teneufm}

\def\boxit#1{\vbox{\hrule\hbox{\vrule\kern2.0pt
       \vbox{\kern2.0pt#1\kern2.0pt}\kern2.0pt\vrule}\hrule}}

\def\vlra#1{\hbox{\kern-1pt
       \hbox{\raise2.38pt\hbox{\vbox{\hrule width#1 height0.26pt}}}
       \kern-4.0pt$\rightarrow$}}

\def\vlla#1{\hbox{$\leftarrow$\kern-1.0pt
       \hbox{\raise2.38pt\hbox{\vbox{\hrule width#1 height0.26pt}}}}}

\def\vlda#1{\hbox{$\leftarrow$\kern-1.0pt
       \hbox{\raise2.38pt\hbox{\vbox{\hrule width#1 height0.26pt}}}
       \kern-4.0pt$\rightarrow$}}

\def\longra#1#2#3{\,\lower3pt\hbox{${\buildrel\mathop{#2}
\over{{\vlra{#1}}\atop{#3}}}$}\,}

\def\longla#1#2#3{\,\lower3pt\hbox{${\buildrel\mathop{#2}
\over{{\vlla{#1}}\atop{#3}}}$}\,}

\def\longda#1#2#3{\,\lower3pt\hbox{${\buildrel\mathop{#2}
\over{{\vlda{#1}}\atop{#3}}}$}\,}

\def\overrightharpoonup#1{\vbox{\ialign{##\crcr
	$\rightharpoonup$\crcr\noalign{\kern-1pt\nointerlineskip}
	$\hfil\displaystyle{#1}\hfil$\crcr}}}

\catcode`@=11
\def\@@dalembert#1#2{\setbox0\hbox{$#1\rm I$}
  \vrule height.90\ht0 depth.1\ht0 width.04\ht0
  \rlap{\vrule height.90\ht0 depth-.86\ht0 width.8\ht0}
  \vrule height0\ht0 depth.1\ht0 width.8\ht0
  \vrule height.9\ht0 depth.1\ht0 width.1\ht0 }
\def\dalembert{\mathord{\mkern2mu\mathpalette\@@dalembert{}\mkern2mu}}

\def\@@varcirc#1#2{\mathord{\lower#1ex\hbox{\m@th${#2\mathchar\hex0017 }$}}}
\def\varcirc{\mathchoice
  {\@@varcirc{.91}\displaystyle}{\@@varcirc{.91}\textstyle}
{\@@varcirc{.45}\scriptscriptstyle}}
\catcode`@=12

\font\tensf=cmss10 \font\sevensf=cmss8 at 7pt
\newfam\sffam
\textfont\sffam=\tensf\scriptfont\sffam=\sevensf

\input amssym.def
\input amssym
\magnification=1200

\font\bigsll=cmsl10 scaled\magstep3
\tolerance=500
\overfullrule=0pt

\null
\bigskip
\centerline{\bigsll Integral models in unramified mixed characteristic (0,2) of }
\bigskip\medskip
\centerline{\bigsll hermitian orthogonal Shimura varieties of PEL type, Part II}
\bigskip\medskip\noindent
\centerline{Adrian Vasiu, Binghamton University}
\bigskip\noindent
\centerline{October 9, 2013, accepted in final form for publication in Math. Nachr.}
\bigskip\noindent
{\bf ABSTRACT.} We construct relative PEL type embeddings in mixed characteristic $(0,2)$ between hermitian orthogonal Shimura varieties of PEL type. We use this to prove the existence of integral canonical models in unramified mixed characteristic $(0,2)$ of hermitian orthogonal Shimura varieties of PEL type. 
\medskip\noindent
{\bf MSC 2010}: 11E57, 11G10, 11G15, 11G18, 14F30, 14G35, 14K10, 14L05, 20G25.
\medskip\noindent
{\bf Key Words}: Shimura varieties, integral models, abelian schemes, $2$-divisible groups, $F$-crystals, reductive and orthogonal group schemes, and involutions.

\bigskip\smallskip
\centerline{\bigsll {\bf 1. Introduction}}
\bigskip\smallskip
This paper is a sequel to [19] and thus its main goal is to prove the existence of {\it integral canonical models} in unramified mixed characteristic $(0,2)$ of {\it hermitian orthogonal Shimura varieties of PEL type}. We begin with a review on Shimura varieties of Hodge in the form needed in the paper. Let ${\bf S}:={\rm Res}_{{\bf C}/{\bf R}} {\bf G}_{m,{\bf C}}$ be the unique two dimensional torus over ${\bf R}$ with the property that ${\bf S}({\bf R})$ is the (multiplicative) group ${\bf G}_{m,{\bf C}}({\bf C})$ of non-zero complex numbers. A Shimura pair $(G,\Mx)$ consists of a reductive group $G$ over ${\bf Q}$ and a $G({\bf R})$-conjugacy class $\Mx$ of homomorphisms ${\bf S}\to G_{\bf R}$ that satisfy Deligne's axioms of [4, Subsubsect. 2.1.1]: (i) the Hodge ${\bf Q}$--structure on the Lie algebra ${\rm Lie}(G)$ of $G$ defined by each element $h\in\Mx$ is of type $\{(-1,1),(0,0),(1,-1)\}$, (ii) no simple factor of the adjoint group $G^{\rm ad}$ of $G$ becomes compact over ${\bf R}$, and (iii) $({\rm Ad}\circ h)({\bf R})(i)$ is a Cartan involution of ${\rm Lie}(G^{\rm ad}_{\bf R})$ in the sense of [8, Ch. III, Sect. 7]. Here ${\rm Ad}:G_{\bf R}\to {\bf GL}_{{\rm Lie}(G^{\rm ad}_{\bf R})}$ is the adjoint representation. The axioms imply that $\Mx$ has a natural structure of a hermitian symmetric domain, cf. [4, Cor. 1.1.17]. 
\smallskip
The most studied Shimura pairs are of the form $({\bf GSp}(W,\psi),\Ms)$, where $(W,\psi)$ is a symplectic space over ${\bf Q}$ and where $\Ms$ is the set of all ${\bf R}$-monomorphisms ${\bf S}\hookrightarrow {\bf GSp}(W,\psi)_{{\bf R}}$ that define Hodge ${\bf Q}$--structures on $W$ of type $\{(-1,0),(0,-1)\}$ and that have either $2\pi i\psi$ or $-2\pi i\psi$ as polarizations. We assume that the Shimura pair $(G,\Mx)$ is of {\it Hodge type} which means that there exists an injective map
$$f:(G,\Mx)\hookrightarrow ({\bf GSp}(W,\psi),\Ms)$$
 of Shimura pairs (see [3], [4], [10, Ch. 1], [11], and [17, Subsect. 2.4]). To $(G,\Mx)$ one associates a complex Shimura variety ${\rm Sh}(G,\Mx)_{\bf C}$ (see [3] and [4]), a number field $E(G,\Mx)$ called the {\it reflex field} (see [3], [4], and [10]), as well as a canonical model ${\rm Sh}(G,\Mx)$ of the complex Shimura variety ${\rm Sh}(G,\Mx)_{\bf C}$ over $E(G,\Mx)$ (see [14] and [3]). 
\smallskip
Let ${\bf Z}_{(2)}$ be the localization of ${\bf Z}$ at the prime ideal $(2)$. Let $L$ be a ${\bf Z}$-lattice of $W$ such that $\psi$ induces a perfect form 
$\psi:L\otimes_{{\bf Z}} L\to{\bf Z}$ i.e., the induced monomorphism $L\hookrightarrow L^\vee:={\rm Hom}(L,{\bf Z})$ is onto. Let $L_{(2)}:=L\otimes_{{\bf Z}} {\bf Z}_{(2)}$. Let $G_{{\bf Z}_{(2)}}$ be the schematic closure of $G$ in the reductive group scheme ${\bf GSp}(L_{(2)},\psi)$ over ${\bf Z}_{(2)}$. Let $K_2:={\bf GSp}(L_{(2)},\psi)({\bf Z}_2)$ and $H_2:=G({\bf Q}_2)\cap K_2=G_{{\bf Z}_{(2)}}({\bf Z}_2)$. Let
$$\Mb:=\{b\in {\rm End}(L_{(2)})|b\;{\rm is}\;{\rm fixed}\; {\rm by}\; G_{{\bf Z}_{(2)}}\}.$$ Let $G_1$ be the subgroup of ${\bf GSp}(W,\psi)$ that fixes all elements of $\Mb[{1\over 2}]$. 
\smallskip
Let $\Mi$ be the involution of the semisimple ${\bf Z}_{(2)}$-algebra ${\rm End}(L_{(2)})$ defined by the identity $\psi(b(l_1),l_2)=\psi(l_1,\Mi(b)(l_2))$,
where $b\in {\rm End}(L_{(2)})$ and $l_1$, $l_2\in L_{(2)}$. As $\Mb=\Mb[{1\over 2}]\cap {\rm End}(L_{(2)})$, we have $\Mi(\Mb)=\Mb$. As the elements of $\Mx$ fix $\Mb\otimes_{\bf Z_{(2)}} {\bf R}$, the involution $\Mi$ of $\Mb$ is positive (in the sense of [9, Sect. 2]). 
\smallskip
Let ${\bf F}$ be an algebraic closure of the field ${\bf F}_2$ with two elements.  Let $W({\bf F})$ be the ring of $2$-typical Witt vectors with coefficients in ${\bf F}$. Let ${\bf A}_f:=\widehat{\bf Z}\otimes_{\bf Z}{\bf Q}$ be the ring of finite ad\`eles. Let ${\bf A}_f^{(2)}$ be the ring of finite ad\`eles with the $2$-component omitted; we have ${\bf A}_f={\bf Q}_2\times{\bf A}_f^{(2)}$. Let $v$ be a prime of $E(G,\Mx)$ that divides $2$. Let $O_{(v)}$ be the localization of the ring of integers of $E(G,\Mx)$ at the prime $v$. 
\medskip\smallskip\noindent
{\bf 1.1. Shimura pairs of PEL type.} In all that follows we will assume that the following four properties (axioms) hold:
\medskip
{\bf (i)} the $W({\bf F})$-algebra $\Mb\otimes_{{\bf Z}_{(2)}} W({\bf F})$ is a product of matrix $W({\bf F})$-algebras;
\smallskip
{\bf (ii)} the ${\bf Q}$--algebra $\Mb[{1\over 2}]$ is ${\bf Q}$--simple;
\smallskip
{\bf (iii)} the group $G$ is the identity component of $G_1$;
\smallskip
{\bf (iv)} the flat, affine group scheme $G_{{\bf Z}_{(2)}}$ over ${\bf Z}_{(2)}$ is reductive (i.e., it is smooth and its special fibre is connected and has a trivial unipotent radical).
\medskip
Property (ii) is not truly required: it is inserted only to ease the presentation. One can check that in fact (iv) implies (i). Property (iv) also implies that $H_2$ is a hyperspecial subgroup of $G({\bf Q}_2)=G_{{\bf Q}_2}({\bf Q}_2)$ (cf. [16, Subsubsect. 3.8.1]) and that the prime $v$ is unramified over $2$ (cf. [11, Cor. 4.7 (a)]). Let $G^{\rm der}$ be the derived group of $G$. In order to avoid trivial cases we will assume that $G$ is  not a torus (equivalently, that $G^{\rm der}$ is non-trivial). Due to properties (ii) and (iii) and this assumption, there exist three possible cases (see [9, Sect. 7]):
\medskip
{\bf (A)} the group $G_{{\bf C}}^{\rm der}$ is a product of ${\bf SL}_n$ groups with $n\ge 2$ and, in the case $n=2$, the center of $G$ has dimension at least $2$;
\smallskip
{\bf (C)} the group $G_{{\bf C}}^{\rm der}$ is a product of ${\bf Sp}_{2n}$ groups with $n\ge 1$ and, in the case $n=1$, the center of $G$ has dimension $1$;
\smallskip
{\bf (D)} the group $G_{{\bf C}}^{\rm der}$ is a product of ${\bf SO}_{2n}$ groups with $n\ge 2$.
\medskip
We have $G\neq G_1$  if and only if we are in the case (D) i.e.,  if and only if $G^{\rm der}$ is not simply connected (cf. [9, Sect. 7]). We recall that PEL stands for polarization, endomorphisms, and level structures and that the first two of these notions refer to the fact that the axiom (iii) holds. In the case (D), one often says that ${\rm Sh}(G,\Mx)$ is a {\it hermitian orthogonal} Shimura variety of PEL type (cf. the description of the intersection group $G_{\bf R}\cap {\bf Sp}(W\otimes_{\bf Q} {\bf R},\psi)$ in [14, Subsects. 2.6 and 2.7]). We are in the case (D)  if and only if $\Mb\otimes_{{\bf Z}_{(2)}} {\bf R}$ is a product of matrix algebras over the quaternion ${\bf R}$-algebra ${\bf H}$ (see [14, Subsect. 2.1, (type III)]). 
\smallskip
In cases (A) and (C), $G_{{\bf Z}_{(2)}}$ is the subgroup scheme of ${\bf GSp}(L_{(2)},\psi)$ that fixes $\Mb$. But in the case (D) one encounters the following problem: the subgroup scheme of ${\bf GSp}(L_{(2)},\psi)$ that fixes $\Mb$ is not smooth and its identity component is not $G_{{\bf Z}_{(2)}}$ (see [19, Subsubsect. 3.5.1 and Remark 3.5.5]). In all that follows we will also assume that:
\medskip
{\bf (v)} we are in the case (D).
\medskip
We refer to each quadruple $(f,L,v,\Mb)$ that satisfies properties (i) to (v) as a {\it hermitian orthogonal standard PEL situation} in mixed characteristic $(0,2)$. 
\smallskip
Let $\Mm$ be the ${\bf Z}_{(2)}$-scheme which parameterizes isomorphism classes of principally polarized abelian schemes that are of relative dimension ${{\dim_{{\bf Q}}(W)}\over 2}$ over ${\bf Z}_{(2)}$-schemes and that are equipped with compatible level-$l$ symplectic similitude structures for all odd numbers $l\in{\bf N}$, cf. [12, Thms. 7.9 and 7.10]. We have a natural identification ${\rm Sh}({\bf GSp}(W,\psi),\Ms)_{E(G,\Mx)}/K_2=\Mm_{E(G,\Mx)}$ as well as an action of ${\bf GSp}(W,\psi)({\bf A}^{(2)}_f)$ on $\Mm$. These symplectic similitude structures and this action are defined naturally via $(L,\psi)$ (see [3, Example 4.16], [11, Sect. 3], and [17, Subsect. 4.1]). 
\smallskip
It is well known that we have an identity ${\rm Sh}(G,\Mx)_{{\bf C}}/H_2=G_{{\bf Z}_{(2)}}({\bf Z}_{(2)})\backslash [\Mx\times G({\bf A}_f^{(2)})]$ and that to the injective map $f$ one associates a natural $E(G,\Mx)$-morphism 
$${\rm Sh}(G,\Mx)/H_2\to {\rm Sh}({\bf GSp}(W,\psi),\Ms)_{E(G,\Mx)}/K_2=\Mm_{E(G,\Mx)}$$ 
which is a closed embedding (for instance, see [19, Subsect. 1.3]). Let $\Mn$ be the schematic closure of ${\rm Sh}(G,\Mx)/H_2$ in $\Mm_{O_{(v)}}$. Let $\Mn^{\rm n}$ be the normalization of $\Mn$. 
\smallskip
Our goal is to prove the following two theorems.
\medskip\smallskip\noindent
{\bf 1.2. Basic Theorem.} {\it Let $(f,L,v,\Mb)$ be a hermitian orthogonal standard PEL situation in mixed characteristic $(0,2)$. Then there exists an embedding $\tilde f:(G,\Mx)\to ({\bf GSp}(\tilde W,\tilde\psi),\tilde\Ms)$ of Shimura pairs and a ${\bf Z}$-lattice $\tilde L$ of $\tilde W$ such that $\tilde f$ factors through an injective map $\tilde f^\prime:(G^\prime,\Mx^\prime)\hookrightarrow ({\bf GSp}(\tilde W,\tilde\psi),\tilde\Ms)$ in such a way that by defining $\tilde L_{(2)}:=\tilde L\otimes_{\bf Z} {\bf Z}_{(2)}$ the following three properties hold:
\medskip
{\bf (i)} the schematic closure of $G$ in ${\bf GL}_{\tilde L_{(2)}}$ is $G_{{\bf Z}_{(2)}}$ and the schematic closure $G^\prime_{{\bf Z}_{(2)}}$ of $G^\prime$ in ${\bf GL}_{\tilde L_{(2)}}$ is a reductive group scheme over ${\bf Z}_{(2)}$ whose extension (pull-back) to ${\bf Z}_2$ is a split ${\bf GSO}_{2nq}$ group scheme for some $q\in {\bf N}$;
\smallskip
{\bf (ii)} if $\tilde\Mb:=\{b\in {\rm End}(\tilde L_{(2)})|b\;{\rm is}\;{\rm fixed}\; {\rm by}\; G_{{\bf Z}_{(2)}}\}$, $v^\prime$ is the prime of $E(G^\prime,\Mx^\prime)$ divided by $v$ (under the natural embedding $E(G^\prime,\Mx^\prime)\hookrightarrow E(G,\Mx)$), and $\tilde\Mb^\prime:=\{b\in {\rm End}(\tilde L_{(2)})|b\;{\rm is}\;{\rm fixed}\; {\rm by}\; G^\prime_{{\bf Z}_{(2)}}\}$, then both quadruples $(\tilde f,\tilde L,v,\tilde\Mb)$ and $(\tilde f^\prime,\tilde L,v^\prime,\tilde\Mb^\prime)$ are hermitian orthogonal standard PEL situations in mixed characteristic $(0,2)$;
\smallskip
{\bf (iii)} the identity component of the subgroup scheme of $G^\prime_{{\bf Z}_{(2)}}$ that fixes $\tilde\Mb$ is $G_{{\bf Z}_{(2)}}$.}
\medskip
We refer to the quadruple $(\tilde f,\tilde f^\prime,\tilde L,v)$ as a {\it hermitian orthogonal relative PEL situation} (it is similar to --though different from-- the {\it relative PEL situations} of [17, Subsubsect. 4.3.16 and Sect. 6]). Property (iii) is the key new ingredient that will allow us to go around the problem mentioned in Subsection 1.1 (i.e., to prove that the $O_{(v)}$-scheme $\Mn^{\rm n}$ is regular and formally smooth) as follows. If $G_{{\bf Z}_2}$ is a split ${\bf GSO}_{2n}$ group scheme, then the $O_{(v)}$-scheme $\Mn$ is regular and formally smooth and thus we have $\Mn=\Mn^{\rm n}$ (cf. [19, Thm. 1.4 (b)]). By combining this result with the Basic Theorem we prove that:
\medskip\smallskip\noindent
{\bf 1.3. Main Theorem.} {\it Let $(f,L,v,\Mb)$ be a hermitian orthogonal standard PEL situation in mixed characteristic $(0,2)$. Then the $O_{(v)}$-scheme $\Mn^{\rm n}$ (obtained as in Subsection 1.1) is regular and formally smooth.} 
\medskip
The locally compact, totally disconnected topological group $G({\bf A}_f^{(2)})$ acts on $\Mn^{\rm n}$ continuously in the sense of [4, Subsubsect. 2.7.1]. Thus $\Mn^{\rm n}$ is an {\it integral canonical model} of ${\rm Sh}(G,\Mx)/H_2$ over $O_{(v)}$ in the sense of [17, Def. 3.2.3 6)], cf. [17, Example 3.2.9 and Cor. 3.4.4]. Due to [20, Cor. 30], $\Mn^{\rm n}$ is the {\it unique} integral canonical model of ${\rm Sh}(G,\Mx)/H_2$ over $O_{(v)}$ and thus as in [17, Rmks. 3.2.4 and 3.2.7 4')] one argues that  $\Mn^{\rm n}$ is the {\it final object} of the category of smooth integral models of ${\rm Sh}(G,\Mx)/H_2$ over $O_{(v)}$; here the word smooth is used as in [10, Def. 2.2]. 
\smallskip
In Section 2 we present tools that pertain to group schemes and that are needed to prove the Basic Theorem in Section 3. In Section 4 we recall few basic crystalline properties from [19]. In Section 5 we prove the Main Theorem. 
\bigskip\smallskip
\centerline{\bigsll {\bf 2. Group schemes}}
\bigskip\smallskip
Let $n\in{\bf N}$. Let ${\rm Spec}\, S$ be an affine scheme. We recall that a reductive group scheme $\Mr$ over $S$ is a smooth, affine group scheme over $S$ whose fibres are connected and have trivial unipotent radicals. Let $\Mr^{\rm ad}$  and $\Mr^{\rm der}$ be the adjoint and the derived (respectively) group schemes of $\Mr$, cf. [6, Vol. III, Exp. XXII, Def. 4.3.6 and Thm. 6.2.1]. Let ${\rm Lie}(\Mu)$ be the Lie algebra of a smooth, closed subgroup scheme $\Mu$ of $\Mr$. For an affine morphism ${\rm Spec}\, S_1\to {\rm Spec}\, S$ and for $Z$ (or $Z_S$ or $Z_*$) an $S$-scheme, let $Z_{S_1}$ (or $Z_{S_1}$ or $Z_{*,S_1}$) be $Z\times_{\rm Spec\, S} {\rm Spec}\, S_1$. If $\tilde S\hookrightarrow S$ is a finite, \'etale ${\bf Z}$-monomorphism, let ${\rm Res}_{S/\tilde S}$ be the operation of Weil restriction from $S$ to $\tilde S$ (see [1, Ch. 7, Sect. 7.6]). Thus ${\rm Res}_{S/\tilde S} \Mr$ is a reductive group scheme over $\tilde S$ such that for each $\tilde S$-scheme $Y$ we have a functorial group identification ${\rm Res}_{S/\tilde S} \Mr(Y)=\Mr(Y\times_{\rm Spec\, \tilde S} {\rm Spec}\, S)$. 
\smallskip
If $M$ is a free $S$-module of finite rank, let $M^\vee:={\rm Hom}_S(M,S)$, let ${\bf GL}_M$ be the reductive group scheme over $S$ of linear automorphisms of $M$, and let $\Mt(M):=\oplus_{s,t\in{\bf N}\cup\{0\}} M^{\otimes s}\otimes_S M^{\vee\otimes t}$. Each $S$-linear isomorphism $\iota:M\tilde\to\tilde M$ of free $S$-modules of finite rank, extends naturally to an $S$-linear isomorphism (to be denoted also by) $\iota:\Mt(M)\tilde\to\Mt(\tilde M)$ and therefore we will speak about $\iota$ taking some tensor of $\Mt(M)$ to some tensor of $\Mt(\tilde M)$. We identify ${\rm End}_S(M)=M\otimes_S M^\vee$. A bilinear form $\lambda_M$ on a free $S$-module $M$ of finite rank is called perfect if it induces an $S$-linear isomorphism $M\tilde\to M^\vee$. If $\lambda_M$ is perfect and alternating, we call the pair $(M,\lambda_M)$ a symplectic space over $S$ and we define ${\bf Sp}(M,\lambda_M):={\bf GSp}(M,\lambda_M)^{\rm der}$. We often use the same notation for two elements of some modules (like involutions, endomorphisms, bilinear forms, etc.) that are obtained one from another via extensions of scalars and restrictions. If $E$ (or $E_*$) is a  number field, let $E_{(2)}$ (or $E_{*,(2)}$) be the normalization of ${\bf Z}_{(2)}$ in $E$ (or $E_*$). 
\smallskip
The reductive group schemes ${\bf Sp}_{2n,S}$, ${\bf GL}_{n,S}$, etc., are over $S$.  The Lie groups ${\bf U}(n)$, ${\bf Sp}(n,{\bf R})$, ${\bf SL}(n,{\bf C})$, ${\bf SO}^*(2n)$, etc., are as in [8, Ch. X, Sect. 2, 1] (but using boldfaced capital letters).  
Let ${\bf SO}_{2n,S}^+$, ${\bf GSO}_{2n,S}^+$, and ${\bf O}_{2n,S}^+$ be the split ${\bf SO}_{2n}$, ${\bf GSO}_{2n}$, and ${\bf O}_{2n}$ (respectively) reductive group schemes over $S$. We recall that ${\bf GSO}_{2n,S}^+$ is the quotient of ${\bf SO}_{2n,S}^+\times_S {\bf G}_{m,S}$ by a $\pmb{\mu}_{2,S}$ subgroup scheme that is embedded diagonally. 
\smallskip
In Subsection 2.1 we review some general facts on ${\bf SO}^+_{2n,{\bf Z}_{(2)}}$. In Subsection 2.2 we study hermitian twists of split ${\bf SO}$ group schemes and the Lie groups associated to them. 
\medskip\smallskip\noindent
{\bf 2.1. Split ${\bf SO}_{2n}$ group schemes.} We consider the quadratic form 
$${\got Q}_n(x):=x_1x_2+\cdots+x_{2n-1}x_{2n}\;\;\; {\rm defined}\;\;\;{\rm for}\;\;\; x=(x_1,\ldots,x_{2n})\in {\got L}_n:={\bf Z}_{(2)}^{2n}.$$ 
For $\alpha\in{\bf Z}_{(2)}$ and $x\in{\got L}_n$ we have ${\got Q}_n(\alpha x)=\alpha^2{\got Q}_n(x)$. Let ${\tilde\Md}_n$ be the subgroup scheme of ${\bf GL}_{{\got L}_n}$ that fixes ${\got Q}_n$. Let $\Md_n$ be the schematic closure of the identity component of ${\tilde{\Md}}_{n,{\bf Q}}$ in $\tilde\Md_n$. We recall from [19, Subsect. 3.1] that $\Md_n$ and $\tilde\Md_n$ are isomorphic to ${\bf SO}^+_{2n,{\bf Z}_{(2)}}$ and ${\bf O}^+_{2n,{\bf Z}_{(2)}}$ (respectively). Thus $\Md_1$ is isomorphic to ${\bf G}_{m,{\bf Z}_{(2)}}$, $\Md_2$ is isomorphic to the quotient of a product of two copies of ${\bf SL}_{2,{\bf Z}_{(2)}}$ by a $\pmb{\mu}_{2,{\bf Z}_{(2)}}$ subgroup scheme that is embedded diagonally, and for $n\ge 2$ the group scheme $\Md_n$ is semisimple. Moreover we have a non-trivial, split short exact sequence $1\to\Md_n\to\tilde\Md_n\to ({{\bf Z}/2{\bf Z}})_{{\bf Z}_{(2)}}\to 1$ and $\Md_n$ is the identity component of $\tilde\Md_n$. Let 
$$\rho_n:\Md_n\hookrightarrow {\bf GL}_{{\got L}_n}$$
be the natural rank $2n$ faithful representation. We also recall from [19, Subsect. 3.1] that $\rho_n$ is associated to the weight $\varpi_1$ if $n\ge 4$ and to the weight $\varpi_2$ of the $A_3$ Lie type if $n=3$ (see [2, plates I and IV] for these weights); moreover $\rho_2$ is the tensor product of the standard rank $2$ representations of the mentioned two copies of ${\bf SL}_{2,{\bf Z}_{(2)}}$. Thus the representation $\rho_n$ is isomorphic to its dual and, up to a ${{\bf G}_m}({\bf Z}_{(2)})$-multiple, there exists a unique perfect, symmetric bilinear form ${\got B}_n$ on ${\got L}_n$ fixed by $\Md_n$ (the case $n=1$ is trivial). In fact we can take ${\got B}_n$ such that
we have ${\got B}_n(u,x):={\got Q}_n(u+x)-{\got Q}_n(u)-{\got Q}_n(x)$ for all $u, x\in {\got L}_n$. Let $J(2n)$ be the matrix representation of ${\got B}_n$ with respect to the standard ${\bf Z}_{(2)}$-basis for ${\got L}_n$; it has $n$ diagonal blocks that are ${0\,1}\choose {1\,0}$. 
\smallskip
Let $q\in {\bf N}$. Let 
$$d_{n,q}:\Md_n^q\hookrightarrow\Md_{nq}$$ 
be a ${\bf Z}_{(2)}$-monomorphism such that the faithful representation $\rho_{nq}\circ d_{n,q}:\Md_n^q\hookrightarrow {\bf GL}_{{\got L}_{nq}}$ is isomorphic to a direct sum of $q$ copies of the representation $\rho_n$; it is easy to see that $d_{n,q}$ is unique up to $\tilde\Md_{nq}({\bf Z}_{(2)})$-conjugation. 
\medskip\noindent
{\bf 2.1.1. Lemma.} {\it Let $r\in {\bf N}$ be a divisor of $q$. Let $\Delta_{n,r,q}:\Md_n^r\hookrightarrow\Md_n^q$ be the diagonal embedding. The composite monomorphism $\rho_{nq}\circ d_{n,q}\circ\Delta_{n,r,q}:\Md_n^r\hookrightarrow\Md_n^q\hookrightarrow \Md_{nq}\hookrightarrow {\bf GL}_{{\got L}_{nq}}$ is a closed embedding that allows us to view $\Md_n^r$ as a closed subgroup scheme of ${\bf GL}_{{\got L}_{nq}}$. Let $\Mb_{n,r,q}$ be the semisimple ${\bf Z}_{(2)}$-subalgebra of ${\rm End}({\got L}_{nq})$ formed by endomorphisms of ${\got L}_{nq}$ fixed by $\Md_n^r$. Then $\Md_n^r$ is the identity component of the centralizer of $\Mb_{n,r,q}$ in $\Md_{nq}$.}
\medskip\noindent
{\it Proof:} The centralizer $\Mc_{n,r,q}$ of $\Mb_{n,r,q}$ in ${\bf GL}_{{\got L}_{nq}}$ is a ${\bf GL}^r_{2n,{\bf Z}_{(2)}}$ group scheme. As $\Mc_{n,r,q}\cap\Md_n^q=\Md_n^r$, to prove the Lemma we can assume that $r=q$. We identify $({\got L}_{nq},{\got B}_{nq})=\oplus_{i=1}^q ({\got L}_n,{\got B}_n)$ in such a way that the identity ${\got L}_{nq}=\oplus_{i=1}^q {\got L}_n$ defines $d_{n,q}$. Then $\Mc_{n,q,q}=\prod_{i=1}^n {\bf GL}_{{\got L}_n}$ and $\Mc_{n,q,q}\cap\tilde\Md_{nq}=\prod_{i=1}^q \tilde\Md_n=\tilde\Md_n^q$. Thus the identity component of the centralizer of $\Mb_{n,q,q}$ in $\Md_{nq}$ is the identity component of $\tilde\Md_n^q$ and therefore it is $\Md_n^q$.\endproof
\medskip\smallskip\noindent
{\bf 2.2. Twists of $\Md_n$.}  In this subsection we list few properties of different twists of $\Md_n$ and of the Lie groups associated to them.
\medskip\noindent
{\bf 2.2.1. On ${\bf SO}^*(2n)$.} Let ${\bf SO}^*(2n)$ be the Lie group over ${\bf R}$ formed by elements of ${\bf SL}(2n,{\bf C})$ that fix the quadratic form $z_1^2+\cdots+z_{2n}^2$ as well as the skew hermitian form $-z_1\bar z_{n+1}+z_{n+1}\bar z_1-\cdots-z_n\bar z_{2n}+z_{2n}\bar z_n$. It is connected (cf. [8, Ch. X, Sect. 2, 2.4]) and it is associated to a reductive (semisimple if $n\ge 2$) group over ${\bf R}$ that is a form of $\Md_{n,{\bf R}}$. 
\smallskip
For $s\in\{1,\ldots,q\}$, let $(z_1^{(s)},\ldots,z_{2n}^{(s)})$ be variables that define an $s$-copy of ${\bf SO}^*(2n)$. The transformation that takes the $2nq$-tuple $(z_1^{(1)},\ldots,z_{2n}^{(1)},\ldots,z_1^{(q)},\ldots,z_{2n}^{(q)})$ to the $2nq$-tuple $(z_1^{(1)},\ldots,z_n^{(1)},\ldots,z_1^{(q)},\ldots,z_n^{(q)},z_{n+1}^{(1)},\ldots,z_{2n}^{(1)},\ldots,z_{n+1}^{(q)},\ldots,z_{2n}^{(q)})$, gives birth to the standard monomorphism of Lie groups 
$$e_{n,q}:{\bf SO}^*(2n)^q\hookrightarrow {\bf SO}^*(2nq).$$ 
By composing the diagonal embedding ${\bf SO}^*(2n)\hookrightarrow {\bf SO}^*(2n)^q$ with $e_{n,q}$ we get a monomorphism:
$$s_{n,q}:{\bf SO}^*(2n)\hookrightarrow {\bf SO}^*(2nq).$$ 
\noindent
{\bf 2.2.2. Lemma.} {\it The following three properties hold:
\medskip
{\bf (a)} we have $s_{n,q}\circ s_{1,n}=s_{1,nq}$;
\smallskip
{\bf (b)} the Lie group ${\bf SO}^*(2)$ is a compact Lie torus of rank $1$;
\smallskip
{\bf (c)} the centralizer $C_{1,n}$ of the image of $s_{1,n}$ in ${\bf SO}^*(2n)$ is a maximal compact Lie subgroup of ${\bf SO}^*(2n)$ isomorphic to ${\bf U}(n)$.} 
\medskip\noindent
{\it Proof:} 
Part (a) is obvious. It is easy to see that ${\bf SO}^*(2)$ is the Lie subgroup of ${\bf SL}(2,{\bf C})$ whose elements take $(z_1,z_2)$ to $(\cos(\theta)z_1+\sin(\theta)z_2,-\sin(\theta)z_1+\cos(\theta)z_2)$ for some $\theta\in{\bf R}$. From this (b) follows. 
\smallskip
We check (c). Let $u_n:{\bf U}(n)\hookrightarrow {\bf SO}^*(2n)$ be the Lie monomorphism that takes $X+iY\in {\bf U}(n)\leqslant {\bf GL}(n,{\bf C})$ to ${{\,\,\,X\,\,\,Y}\choose {-Y\,X}}\in {\bf SO}^*(2n)\leqslant {\bf GL}(2n,{\bf C})$, where both $X$ and $Y$ are real $n\times n$ matrices (see [8, Ch. X, Sect. 2, 3, Type D III]). The image through $u_n$ of the center of ${\bf U}(n)$ is ${\rm Im}(s_{1,n})$. Thus we have ${\rm Im}(u_n)\leqslant C_{1,n}$. But the centralizer of the image of the complexification of $s_{1,n}$ in the complexification of ${\bf SO}^*(2n)$ is a ${\bf GL}(n,{\bf C})$ Lie group, cf. the very definition of $s_{1,n}$. Thus $C_{1,n}$  is a form of ${\bf GL}(n,{\bf R})$. By reasons of dimensions we get that $C_{1,n}={\rm Im}(u_n)\tilde\to {\bf U}(n)$. The fact that $C_{1,n}$ is a maximal compact Lie subgroup of ${\bf SO}^*(2n)$ follows from [8, Ch. X, Sect. 6, Table V]. Thus (c) holds.\endproof
\medskip\noindent
{\bf 2.2.3. Hermitian twists.} Let ${\rm Spec}\, S_2\to {\rm Spec}\, S_1$ be an \'etale cover of degree $2$ between regular schemes that are flat over ${\bf Z}_{(2)}$ and with $S_1$ integral. There exists a unique $S_1$-automorphism $\tau\in {\rm Aut}_{S_1}(S_2)$ of order $2$. Let $M_1:=S_1^{2n}$ and $M_2:=S_2^{2n}=M_1\otimes_{S_1} S_2$. We view also ${\got Q}_n(x)$ as a quadratic form defined for $x=(x_1,\ldots,x_{2n})\in M_2$. Let
$${\got H}_n^*(x):=x_1\tau(x_1)-x_2\tau(x_2)+\cdots+x_{2n-1}\tau(x_{2n-1})-x_{2n}\tau(x_{2n});$$ 
it is a skew hermitian form with respect to $\tau$ on $M_2=S_2^{2n}$. Let $\Md_{n,\tau}$ be the subgroup scheme of ${\rm Res}_{S_2/S_1} \Md_{n,S_2}$ that fixes ${\got H}_n^*$. 
\medskip\noindent
{\bf 2.2.4. Lemma.} {\it We have the following four properties:
\medskip
{\bf (a)} the group scheme $\Md_{n,\tau}$ over $S_1$ is reductive, splits over $S_2$, and for $n\ge 2$ it is semisimple;
\smallskip
{\bf (b)} if $S_2=S_1\oplus S_1$ as $S_1$-algebras, then $\Md_{n,\tau}$ is isomorphic to $\Md_{n,S_1}$ and thus it is split;
\smallskip
{\bf (c)} if $S_1\hookrightarrow S_2$ is ${\bf R}\hookrightarrow{\bf C}$ (thus $\tau$ is the complex conjugation), then $\Md_{n,\tau}({\bf R})$ is isomorphic to ${\bf SO}^*(2n)$;
\smallskip
{\bf (d)} there exist $S_1$-monomorphisms $\Md_{n,\tau}\hookrightarrow {\bf Sp}_{8n,S_1}$ such that the resulting rank $8n$ representation of $\Md_{n,\tau,S_2}$ is isomorphic to four copies of the representation $\rho_{n,S_2}$.} 
\medskip\noindent
{\it Proof:}
To prove (a) and (b) we can assume that $S_2=S_1\oplus S_1$ as $S_1$-algebras; thus $\tau$ is the permutation of the two factors $S_1$ of $S_2$. For $j\in\{1,\ldots,2n\}$ we write 
$$x_j=(u_j,(-1)^{j+1}v_j)\in S_2=S_1\oplus S_1$$ 
with $u_j,v_j\in S_1$. Thus we have identities
${\got Q}_n(x)=(u_1u_2+\cdots+u_{2n-1}u_{2n},v_1v_2+\cdots+v_{2n-1}v_{2n})$ and ${\got H}_n^*(x)=(\sum_{j=1}^{2n} u_jv_j,\sum_{j=1}^{2n} u_jv_j)$. Let 
$$(g_1,g_2)\in {\rm Res}_{S_2/S_1} \Md_{n,S_2}(S_1)=\Md_n(S_1)\times \Md_n(S_1)\leqslant {\bf GL}_{M_1}(S_1)\times {\bf GL}_{M_1}(S_1)$$ 
be such that it fixes $\sum_{j=1}^{2n} u_jv_j$; here $g_1$ and $g_2$ act on an $M_1$ copy that involves the variables $(u_1,\ldots,u_{2n})$ and $(v_1,\ldots,v_{2n})$ (respectively). We get that at the level of $2n\times 2n$ matrices we have $g_2=(g_1^t)^{-1}=J(2n)g_1J(2n)$. Thus the $S_1$-monomorphism $\Md_{n,\tau}\hookrightarrow {\rm Res}_{S_2/S_1} \Md_{n,S_2}$ is isomorphic to the diagonal $S_1$-monomorphism $\Md_{n,S_1}\hookrightarrow \Md_{n,S_1}\times_{S_1} \Md_{n,S_1}$ and therefore $\Md_{n,\tau}$ is isomorphic to $\Md_{n,S_1}$. From this (a) and (b) follow. 
\smallskip
We check (c). Let $i\in{\bf C}$ be the standard square root of $-1$ and let $\bar w:=\tau(w)$, where $w\in{\bf C}$. Under the transformation $x_{2j-1}:=z_j+iz_{n+j}$ and $x_{2j}:=z_j-iz_{n+j}$ (with $j\in\{1,\ldots,n\}$), we have 
$${\got Q}_n(x)=z_1^2+z_2^2+\cdots+z_{2n}^2\;\; {\rm and}\;\; {\got H}_n^*(x)=2i(-z_1\bar z_{n+1}+z_{n+1}\bar z_1-\cdots-z_n\bar z_{2n}+z_{2n}\bar z_n).$$ 
From the very definition of ${\bf SO}^*(2n)$ we get that (c) holds. 
\smallskip
For $s\in {\bf N}$ we have standard closed embedding monomorphisms 
$${\bf GL}_{sn,S_2}\hookrightarrow {\bf Sp}_{2sn,S_2}\;\; {\rm and}\;\; {\rm Res}_{S_2/S_1} {\bf Sp}_{2sn,S_2}\hookrightarrow {\bf Sp}_{4sn,S_1}.$$ Therefore ${\rm Res}_{S_2/S_1} \Md_{n,S_2}$ is (via ${\rm Res}_{S_2/S_1} \rho_{n,S_2}$) naturally a closed subgroup scheme of ${\rm Res}_{S_2/S_1} {\bf GL}_{2n,S_2}$ and therefore also of ${\rm Res}_{S_2/S_1} {\bf Sp}_{4n,S_2}$ and of ${\bf Sp}_{8n,S_1}$. Thus (d) holds.\endproof
\medskip\noindent
{\bf 2.2.5. The case of number fields.} Let $F_0$ be a totally imaginary quadratic extension of ${\bf Q}$ in which $2$ splits. Let $S_1:={\bf Z}_{(2)}$ and $S_2:=F_{0,(2)}$; thus $\tau$ is the non-trivial ${\bf Z}_{(2)}$-automorphism of $F_{0,(2)}$ and the reductive group scheme $\Md_{n,\tau}$ is a form of $\Md_n$. As $2$ splits in $F_0$ and as $\Md_{n,\tau}$ splits over $F_{0,(2)}$ (cf. Lemma 2.2.4 (b)), $\Md_{n,\tau}$ splits over ${\bf Z}_2$. Let $F_1$ be a totally real, finite Galois extension of ${\bf Q}$ unramified above $2$ and of degree $q\in {\bf N}$. Let $F_2:=F_1\otimes_{{\bf Q}} F_0$; it is a totally imaginary, Galois extension of ${\bf Q}$ that has degree $2q$ and that is unramified above $2$. We view $N_1:=F_{1,(2)}^{2n}$ as a free ${\bf Z}_{(2)}$-module of rank $2nq$ and $N_2:=F_{2,(2)}^{2n}=N_1\otimes_{\bf Z_{(2)}} F_{0,(2)}$ as a free $F_{0,(2)}$-module of rank $2nq$. 
\smallskip
See Subsection 2.1 and Subsubsection 2.2.3 for ${\got Q}_n$ and ${\got H}_n^*$ (respectively). Let ${\got Q}_n^{F_{1,(2)}/{\bf Z}_{(2)}}$ and ${\got H}_n^{*F_{1,(2)}/{\bf Z}_{(2)}}$ be (the tensorizations with $F_{1,(2)}$ over ${\bf Z}_{(2)}$ of) ${\got Q}_n$ and ${\got H}_n^*$ (respectively) but viewed as a quadratic form in $2nq$ variables on $N_2$ and as a skew hermitian form with respect to $\tau$ in $2nq$ variables on $N_2$ (respectively). If we fix an $F_{0,(2)}$-linear isomorphism $c_2:N_2\tilde\to F_{0,(2)}^{2nq}$, then for $w=(w_1,\ldots,w_{2nq})\in N_2\tilde\to F_{0,(2)}^{2nq}$ we have 
$${\got Q}_n^{F_{1,(2)}/{\bf Z}_{(2)}}(w)={\got Q}_n(x)\;\; {\rm and}\;\; {\got H}_n^{*F_{1,(2)}/{\bf Z}_{(2)}}(w)={\got H}_n^*(x),$$ 
where $x=(x_1,\ldots,x_{2n}):=w\in F_{2,(2)}^{2n}=N_2$ is computed via the standard $F_{2,(2)}$-basis for $F_{2,(2)}^{2n}$. By taking $c_2$ to be the natural tensorization over $F_{0,(2)}$ of a ${\bf Z}_{(2)}$-linear isomorphism $c_1:N_1\tilde\to{\bf Z}_{(2)}^{2nq}$, we can also view naturally ${\got Q}_n^{F_{1,(2)}/{\bf Z}_{(2)}}$ as a quadratic form in $2nq$ variables on the ${\bf Z}_{(2)}$-module $N_1$. 
\smallskip
Let $\Md^{F_{1,(2)}/{\bf Z}_{(2)}}_{nq,\tau}$ be the identity component of the subgroup scheme of ${\rm Res}_{F_{0,(2)}/{\bf Z}_{(2)}} {\bf GL}_{N_2}$ that fixes both ${\got Q}_n^{F_{1,(2)}/{\bf Z}_{(2)}}$ and ${\got H}_n^{*F_{1,(2)}/{\bf Z}_{(2)}}$; it is a group scheme over ${\bf Z}_{(2)}$.  
\medskip\noindent
{\bf 2.2.6. Lemma.} {\it The following four properties hold:
\medskip
{\bf (a)} the group scheme $(\Md^{F_{1,(2)}/{\bf Z}_{(2)}}_{nq,\tau})_{F_{1,(2)}}$ is isomorphic to $\Md_{nq,\tau,F_{1,(2)}}$; 
\smallskip
{\bf (b)} the group scheme $\Md^{F_{1,(2)}/{\bf Z}_{(2)}}_{nq,\tau}$ is reductive for $n\ge 1$ and semisimple for $n\ge 2$;
\smallskip
{\bf (c)} the group scheme $\Md^{F_{1,(2)}/{\bf Z}_{(2)}}_{nq,\tau}$ splits over $F_{0,(2)}$ and thus also over ${\bf Z}_2$;
\smallskip
{\bf (d)} if $(\tilde L_{(2)},\tilde\psi)$ is  a symplectic space over ${\bf Z}_{(2)}$ of rank $8nq$, then there exist ${\bf Z}_{(2)}$-monomorphisms of reductive group schemes
$${\rm Res}_{F_{1,(2)}/{\bf Z}_{(2)}} \Md_{n,\tau,F_{1,(2)}}\hookrightarrow \Md_{nq,\tau}^{F_{1,(2)}/{\bf Z}_{(2)}}\hookrightarrow {\bf Sp}(\tilde L_{(2)},\tilde\psi)\leqno (1)$$
which are closed embeddings and for which the following three things hold:
\medskip
\item{\bf (i)} the Lie monomorphism ${\rm Res}_{F_{1,(2)}/{\bf Z}_{(2)}} \Md_{n,\tau,F_{1,(2)}}({\bf R})\hookrightarrow \Md_{nq,\tau}^{F_{1,(2)}/{\bf Z}_{(2)}}({\bf R})$ is isomorphic to the Lie monomorphism $e_{n,q}:{\bf SO}^*(2n)^q\hookrightarrow {\bf SO}^*(2nq)$ of Subsubsection 2.2.1;
\smallskip
\item{\bf (ii)} the extension of ${\rm Res}_{F_{1,(2)}/{\bf Z}_{(2)}} \Md_{n,\tau,F_{1,(2)}}\hookrightarrow \Md_{nq,\tau}^{F_{1,(2)}/{\bf Z}_{(2)}}$ to $F_{2,(2)}$ is isomorphic to the extension of the monomorphism $d_{n,q}:\Md_n^q\hookrightarrow \Md_{nq}$ (of Subsection 2.1) to $F_{2,(2)}$;
\smallskip
\item{\bf (iii)} the faithful representation $(\Md_{nq,\tau}^{F_{1,(2)}/{\bf Z}_{(2)}})_{F_{0,(2)}}\hookrightarrow {\bf GL}_{\tilde L_{(2)}\otimes_{{\bf Z}_{(2)}} F_{0,(2)}}$ is isomorphic to four  copies of the representation $\rho_{nq,F_{0,(2)}}$. }
\medskip\noindent
{\it Proof:}
As $F_1$ and $F_2$ are Galois extensions of ${\bf Q}$ unramified above $2$, we have natural identifications $F_{1,(2)}\otimes_{{\bf Z}_{(2)}} F_{1,(2)}=F_{1,(2)}^q$ and $F_{2,(2)}\otimes_{{\bf Z}_{(2)}} F_{1,(2)}=F_{2,(2)}^q$ of ${\bf Z}_{(2)}$-algebras. Thus for 
$$x=(x_1,\ldots,x_{2n})\in N_2\otimes_{{\bf Z}_{(2)}} F_{1,(2)}=(F_{2,(2)}\otimes_{{\bf Z}_{(2)}} F_{1,(2)})^{2n}$$ 
we can speak about the transformation 
$$x_{2j-\epsilon}=(w_{(2j-2)q+2-\epsilon},w_{(2j-2)q+4-\epsilon},\ldots,w_{(2j-2)q+2q-\epsilon})\in F_{2,(2)}\otimes_{{\bf Z}_{(2)}} F_{1,(2)}=F_{2,(2)}^q,$$ 
where $j\in\{1,\ldots,n\}$, $\epsilon\in\{0,1\}$, and $w:=(w_1,\ldots,w_{2nq})\in F_{2,(2)}^{2nq}$. Thus by considering the composite $F_{1,(2)}$-linear isomorphism $N_2\otimes_{{\bf Z}_{(2)}} F_{1,(2)}=F_{0,(2)}^{2nq}\otimes_{{\bf Z}_{(2)}} F_{1,(2)}\tilde\to F_{1,(2)}\otimes_{{\bf Z}_{(2)}} F_{0,(2)}^{2nq}=F_{2,(2)}^{2nq}$ whose inverse is defined naturally by this transformation, we have the following identities 
$${\got Q}_n^{F_{1,(2)}/{\bf Z}_{(2)}}(w)=\sum_{s=1}^{nq} w_{2s-1}w_{2s}\;\;\;{\rm and}\;\;\;\;{\got H}_n^{*F_{1,(2)}/{\bf Z}_{(2)}}(w)=\sum_{s=1}^{2nq} (-1)^{s+1}w_s(1_{F_{1,(2)}}\otimes \tau)(w_s).$$ 
Thus we can redefine $(\Md^{F_{1,(2)}/{\bf Z}_{(2)}}_{nq,\tau})_{F_{1,(2)}}$ as the subgroup scheme of ${\rm Res}_{F_{2,(2)}/F_{1,(2)}} \Md_{nq,F_{2,(2)}}$ that fixes ${\got H}_{nq}^*$. From this (a) follows. Part (b) is implied by (a) and Lemma 2.2.4 (a). 
\smallskip
To check (c), we first remark that we have an identification $F_{2,(2)}\otimes_{{\bf Z}_{(2)}} F_{0,(2)}=F_{2,(2)} \oplus F_{2,(2)}$ of ${\bf Z}_{(2)}$-algebras. Thus, similar to the proof of Lemma 2.2.4 (b) we get that $(\Md^{F_{1,(2)}/{\bf Z}_{(2)}}_{nq,\tau})_{F_{0,(2)}}$ is isomorphic to the identity component of the subgroup scheme of ${\bf GL}_{N_2}={\bf GL}_{N_1\otimes_{{\bf Z}_{(2)}} F_{0,(2)}}$ that fixes ${\got Q}_n^{F_{1,(2)}/{\bf Z}_{(2)}}$ but viewed as a quadratic form in $2nq$ variables on the $F_{0,(2)}$-module $N_2=N_1\otimes_{{\bf Z}_{(2)}} F_{0,(2)}$. We choose an $F_{0,(2)}$-basis $\{y_1,y_2,,\ldots,y_{2nq}\}$ for $N_2=N_1\otimes_{{\bf Z}_{(2)}} F_{0,(2)}$ such that for $\epsilon\in\{0,1\}$ the set $\{y_{2-\epsilon},y_{4-\epsilon},\ldots,y_{2nq-\epsilon}\}$ is formed by elements of the $2-\epsilon$-th, $\ldots$, $2q-\epsilon$-th $F_{1,(2)}\otimes_{{\bf Z}_{(2)}} F_{0,(2)}$ copy of $N_2=N_1\otimes_{{\bf Z}_{(2)}} F_{0,(2)}=(F_{1,(2)}\otimes_{{\bf Z}_{(2)}} F_{0,(2)})^{2n}$. If $i_1$, $i_2\in\{1,\ldots,2nq\}$ are congruent modulo $2$, then we have 
$${\got B}_n^{F_{1,(2)}/{\bf Z}_{(2)}}(y_{i_1},y_{i_2}):={\got Q}_n^{F_{1,(2)}/{\bf Z}_{(2)}}(y_{i_1}+y_{i_2})-{\got Q}_n^{F_{1,(2)}/{\bf Z}_{(2)}}(y_{i_1})-{\got Q}_n^{F_{1,(2)}/{\bf Z}_{(2)}}(y_{i_2})=0.$$ 
This implies the existence of an $F_{0,(2)}$-basis $\{y_1^\prime,y_2,\ldots,y_{2nq-1}^\prime,y_{2nq}\}$ for $N_2=N_1\otimes_{{\bf Z}_{(2)}} F_{0,(2)}$ that has the following two properties:
\medskip
{\bf (iv)} the $F_{0,(2)}$-spans of $\{y_1,y_3,\ldots,y_{2q-1}\}$ and $\{y_1^\prime,y_3^\prime,\ldots,y_{2q-1}^\prime\}$ are equal;
\smallskip
{\bf (v)} the matrix representation of the bilinear form ${\got B}_n^{F_{1,(2)}/{\bf Z}_{(2)}}$ defined by ${\got Q}_n^{F_{1,(2)}/{\bf Z}_{(2)}}$ with respect to the $F_{0,(2)}$-basis $\{y_1^\prime,y_2,\ldots,y_{2nq-1}^\prime,y_{2nq}\}$ for $N_2=N_1\otimes_{{\bf Z}_{(2)}} F_{0,(2)}$, is $J(2nq)$.
\medskip
From (v) we get that $(\Md^{F_{1,(2)}/{\bf Z}_{(2)}}_{nq,\tau})_{F_{0,(2)}}$ is isomorphic to $\Md_{nq,F_{0,(2)}}$ and thus it is a split group scheme. Therefore (c) holds.
\smallskip
We have a canonical ${\bf Z}_{(2)}$-monomorphism ${\rm Res}_{F_{1,(2)}/{\bf Z}_{(2)}} \Md_{n,\tau,F_{1,(2)}}\hookrightarrow \Md_{nq,\tau}^{F_{1,(2)}/{\bf Z}_{(2)}}$; more precisely, the group scheme ${\rm Res}_{F_{1,(2)}/{\bf Z}_{(2)}} \Md_{n,\tau,F_{1,(2)}}$ is the closed subgroup scheme of $\Md_{nq,\tau}^{F_{1,(2)}/{\bf Z}_{(2)}}$ whose group of ${\bf Z}_{(2)}$-valued points is the maximal subgroup of $\Md_{nq,\tau}^{F_{1,(2)}/{\bf Z}_{(2)}}({\bf Z}_{(2)})$ which is formed by $F_{2,(2)}$-linear automorphisms of $N_2=F_{2,(2)}^{2n}$. We take $\Md_{nq,\tau}^{F_{1,(2)}/{\bf Z}_{(2)}}\hookrightarrow {\bf Sp}(\tilde L_{(2)},\tilde\psi)$ to be the composite of the following five ${\bf Z}_{(2)}$-monomorphisms 
$$\Md_{nq,\tau}^{F_{1,(2)}/{\bf Z}_{(2)}}\hookrightarrow {\rm Res}_{F_{0,(2)}/{\bf Z}_{(2)}} (\Md_{nq,\tau}^{F_{1,(2)}/{\bf Z}_{(2)}})_{F_{0,(2)}}$$
$$\tilde\to {\rm Res}_{F_{0,(2)}/{\bf Z}_{(2)}} \Md_{nq,F_{0,(2)}}\hookrightarrow {\rm Res}_{F_{0,(2)}/{\bf Z}_{(2)}} {\bf Sp}_{4nq,F_{0,(2)}}\hookrightarrow {\bf Sp}_{8nq,{\bf Z}_{(2)}}\tilde\to {\bf Sp}(\tilde L_{(2)},\tilde\psi)$$ 
(cf. (c) for the second one and proof of Lemma 2.2.4 (d) for the third and fourth ones). 
\smallskip
Part (i) of (d) follows from the proof of (a): the variables $z_1^{(1)}$, $\ldots$, $z_{2n}^{(q)}$ we used to define the Lie monomorphism $e_{n,q}$ are the variables $w_1$, $\ldots$, $w_{2nq}$ up to a natural permutation. The extension of ${\rm Res}_{F_{1,(2)}/{\bf Z}_{(2)}} \Md_{n,\tau,F_{1,(2)}}$ to $F_{2,(2)}$ is isomorphic to ${\rm Res}_{F_{2,(2)}^q/F_{2,(2)}} \Md_{n,\tau,F_{2,(2)}}\tilde\to \Md_{n,\tau,F_{2,(2)}}^q$ and therefore also to $\Md_{n,F_{2,(2)}}^q$, cf. Lemma 2.2.4 (b). Moreover, the extension of $\Md_{nq,\tau}^{F_{1,(2)}/{\bf Z}_{(2)}}$ to $F_{2,(2)}$ is isomorphic to $\Md_{nq,F_{2,(2)}}$ (cf. (a) and (c)). Thus part (ii) of (d) follows easily from constructions. Part (iii) of (d) follows from (c) and Lemma 2.2.4 (d), cf. constructions.\endproof
\bigskip\smallskip
\centerline{\bigsll {\bf 3. Proof of the Basic Theorem}}
\bigskip\smallskip
In this section the following list of notations 
$${\bf S},\,f:(G,\Mx)\hookrightarrow ({\bf GSp}(W,\psi),\Ms),\,E(G,\Mx),v,\,O_{(v)},\, L,\, L_{(2)},\,G_{\bf Z_{(2)}},\,K_2,\,H_2,\,\Mb,\,\Mi\leqno ({\sharp}_1)$$ 
will be as in Section 1. In order to prove the Basic Theorem (see Subsections 3.3 to 3.5), we will need few extra notations and properties (see Lemma 3.1 and Subsection 3.2). In Subsection 3.3 we use a twisting process to construct using Lemma 2.2.6 the closed embedding monomorphisms $G_{{\bf Z}_{(2)}}\hookrightarrow G_{{\bf Z}_{(2)}}^{\prime}\hookrightarrow {\bf GSp}(\tilde L_{(2)},\tilde\psi)$. In Subsection 3.4 we show that indeed we have injective maps $(G,\Mx)\hookrightarrow (G^\prime,\Mx^\prime)\hookrightarrow ({\bf GSp}(\tilde W,\tilde\psi)$ between Shimura pairs. The proof of the Basic Theorem 1.2 is completed in Subsection 3.5. 
\smallskip
Let $n$, $r\in {\bf N}$ with $n\ge 2$ be such that the group $G_{{\bf C}}^{\rm der}$ is isomorphic to ${\bf SO}^r_{2n,{\bf C}}$, cf. property 1.1 (v). Let $d:={{rn(n-1)}\over 2}\in {\bf N}$. Let $B({\bf F})$ be the field of fractions of $W({\bf F})$. 
\smallskip
Let $\Mb_1$ be the centralizer of $\Mb$ in ${\rm End}(L_{(2)})$. Let $G_{2,{\bf Z}_{(2)}}$ be the centralizer of $\Mb$ in ${\bf GL}_{L_{(2)}}$; thus $G_{2,{\bf Z}_{(2)}}$ is the reductive group scheme over ${\bf Z}_{(2)}$ of invertible elements of $\Mb_1$ and ${\rm Lie}(G_{2,{\bf Z}_{(2)}})$ is the Lie algebra associated to $\Mb_1$. Due to the property 1.1 (i), the $W({\bf F})$-algebra $\Mb_1\otimes_{{\bf Z}_{(2)}} W({\bf F})$ is also a product of matrix $W({\bf F})$-algebras. This implies that $G_{2,W({\bf F})}$ is a product of ${\bf GL}_{2n}$ group schemes over $W({\bf F})$. For basic terminology on involutions of semisimple algebras we refer to [19, Subsect. 3.3]. All simple factors of $(\Mb_j\otimes_{{\bf Z}_2} W({\bf F}),\Mi)$ have the same type which (due to the property 1.1 (v)) is the first orthogonal type. Thus the $W({\bf F})$-monomorphism $G^{\rm der}_{W({\bf F})}\hookrightarrow G_{2,W({\bf F})}$ is isomorphic to the product of $r$ copies of $\rho_{n,W({\bf F})}$ and thus $G_{2,W({\bf F})}$ is isomorphic to ${\bf GL}^r_{2n,W({\bf F})}$. 
\medskip\noindent\noindent
{\bf 3.1. Lemma.} {\it There exists a totally real number field $F$ of degree $r$ and a semisimple group $\Mg$ over $F$ such that the following two properties hold:
\medskip
{\bf (i)} the derived group $G^{\rm der}$ is isomorphic to ${\rm Res}_{F/{\bf Q}} \Mg$;
\smallskip
{\bf (ii)} the group $\Mg$ is a form of ${\bf SO}^+_{2n,F}$ (i.e., of $\Md_{n,F}$).} 
\medskip\noindent
{\it Proof:} Let $F$ be the center of $\Mb_1[{1\over 2}]$. Due to the axiom 1.1 (ii), the ${\bf Q}$--algebra $\Mb_1[{1\over 2}]$ is simple. Thus $F$ is a number field. As $G_{2,W({\bf F})}$ is isomorphic to ${\bf GL}^r_{2n,W({\bf F})}$, we have $[F:{\bf Q}]=r$. Let $G_3$ be the reductive group over $F$ of invertible elements of (warning!) the $F$-algebra $\Mb_1[{1\over 2}]$. We can identify $G_2$ with ${\rm Res}_{F/{\bf Q}} G_3$. As $G^{\rm der}$ is a subgroup of $G_2$ such that $G^{\rm der}_{\bf C}$ is a product of subgroups of the factors of the following product $G_{2,\bf C}=({\rm Res}_{F/{\bf Q}} G_3)_{\bf C}=\prod_{j:F\hookrightarrow {\bf C}} G_3\times_F {}_j {\bf C}$, there exists a subgroup $\Mg$ of $G_3$ such that we can identify $G^{\rm der}$ with ${\rm Res}_{F/{\bf Q}} \Mg$. More precisely, we have a natural product decomposition $F\otimes_{\bf Q} F=F\times F^\perp$ of \'etale $F$-algebras and therefore we can identify $G_{2,F}$ with $G_3\times_F {\rm Res}_{F^\perp/F} G_{3,F^\perp}$; we can now take $\Mg:=G_3\cap G^{\rm der}_F$, the intersection being taken inside $G_{2,F}$. Thus (i) holds. 
\smallskip
It is easy to see that each group $\Mg\times_F {}_j {\bf C}$ is isomorphic to ${\bf SO}_{2n,{\bf C}}$. Thus (ii) also holds. As each simple factor of $G^{\rm ad}_{\bf R}\tilde\to {\rm Res}_{F\otimes_{\bf Q} {\bf R}/{\bf R}} (\Mg\times_F (F\otimes_{\bf Q} {\bf R}))$ is absolutely simple (cf. [4, Subsubsect. 2.3.4 (a)]), the ${\bf R}$-algebra $F\otimes_{\bf Q} {\bf R}$ is isomorphic to ${\bf R}^r$. Thus the number field $F$ is indeed totally real.\endproof 
\medskip\smallskip\noindent
{\bf 3.2. Constructing $G_{{\bf Z}_{(2)}}\hookrightarrow G_{{\bf Z}_{(2)}}^{\prime}\hookrightarrow {\bf GSp}(\tilde L_{(2)},\tilde\psi)$.}  Let $\kappa$ be the set of primes of $F$ that divide $2$. We have a  natural product decomposition 
$$F\otimes_{{\bf Q}} {\bf Q}_2=\prod_{j\in \kappa} F_j$$ 
into $2$-adic fields. Due to Lemma 3.1 (i), we can identify $G_{{\bf Q}_2}^{\rm der}$ with $\prod_{j\in \kappa} {\rm Res}_{F_j/{\bf Q}_2} \Mg_{F_j}$. Thus we can also identify $G^{\rm der}_{B({\bf F})}$ with $\prod_{j\in \kappa} {\rm Res}_{F_j\otimes_{{\bf Q}_2} B({\bf F})/B({\bf F})} (\Mg_{F_j}\times_{F_j} (F_j\otimes_{{\bf Q}_2} B({\bf F})))$. As $G_{\bf Z_{(2)}}$ splits over $W({\bf F})$, the group $G^{\rm der}_{B({\bf F})}$ is split. Thus the $B({\bf F})$-algebra $F_j\otimes_{{\bf Q}_2} B({\bf F})$ is isomorphic to a product of copies of $B({\bf F})$. Therefore each field $F_j$ is unramified over ${\bf Q}_2$ i.e., $F$ is unramified above $2$. Thus the finite ${\bf Z_{(2)}}$-algebra $F_{(2)}$ (see Section 2) is \'etale. 
\smallskip
As $\Mb_1$ is a semisimple ${\bf Z}_{(2)}$-algebra, it is also a semisimple $F_{(2)}$-algebra. Let $G_{3,F_{(2)}}$ be the reductive group scheme over $F_{(2)}$ of invertible elements of the $F_{(2)}$-algebra $\Mb_1$. We can identify $G_{2,{\bf Z}_{(2)}}$ with ${\rm Res}_{F_{(2)}/{\bf Z}_{(2)}} G_{3,F_{(2)}}$. Let $\Mg_{F_{(2)}}$ be the schematic closure of $\Mg$ in $G_{3,F_{(2)}}$. We can identify $G^{\rm der}_{{\bf Z}_{(2)}}$ with ${\rm Res}_{F_{(2)}/{\bf Z}_{(2)}} \Mg_{F_{(2)}}$ and this implies that $\Mg_{F_{(2)}}$ is a semisimple group scheme over $F_{(2)}$. 
\smallskip
Let $O_j$ be the ring of integers of the $2$-adic field $F_j$. For a later use we point out that we can identify $G^{\rm der}_{{\bf Z}_2}$ with $\prod_{j\in\kappa} {\rm Res}_{O_j/{\bf Z}_2} \Mg_{O_j}$. Thus we can also identify $\kappa$ with the set of factors of $G^{\rm der}_{{\bf Z}_2}$ that are Weil restrictions of semisimple ${\bf SO}_{2n}$ group schemes. 
\medskip\smallskip\noindent
{\bf 3.3. A twisting process.}  We fix a totally imaginary quadratic extension $F_0$ of ${\bf Q}$ in which $2$ splits. Let $\tau$ be the non-trivial automorphism of $F_{0,(2)}$ and let $\Md_{n,\tau}$ be the form of $\Md_n$ introduced in Subsubsection 2.2.5. The Lie group $G^{\rm der}_{\bf R}({\bf R})$ is isomorphic to ${\bf SO}^*(2n)^r$, cf. [14, Subsects. 2.6 and 2.7]. This property implies that each connected component of $\Mx$ is a product of $r$ copies of the irreducible hermitian symmetric domain associated to ${\bf SO}^*(2n)$ and thus (cf. [8, Ch. X, Sect. 6, Table V]) we have 
$$\dim_{\bf C}(\Mx)=d=rn(n-1)/ 2.$$ 
The mentioned property also implies that the semisimple group schemes $\Mg_{F_{(2)}}$ and $\Md_{n,\tau,F_{(2)}}$ become isomorphic under extensions via ${\bf Z}_{(2)}$-monomorphisms $F_{(2)}\hookrightarrow{\bf R}$, cf. also the property (i) of Lemma 2.2.6. Thus the class $\gamma\in H^1(F,{\it Aut}(\Mg^{\rm ad}_F))$ that defines the form $\Md_{n,\tau,F}$ of $\Mg$ maps into the trivial class of $H^1({\bf R},{\it Aut}(\Mg^{\rm ad}_{\bf R}))$ via all embeddings $F\hookrightarrow {\bf R}$.
\smallskip
Let ${\it Aut}^\prime(\Mg^{\rm ad}_{F_{(2)}})$ be the smooth, closed, normal subgroup scheme of the group scheme of automorphisms ${\it Aut}(\Mg^{\rm ad}_{F_{(2)}})$ such that we have a first short exact sequence 
$$1\to {\it Aut}^\prime (\Mg^{\rm ad}_{F_{(2)}})\to {\it Aut} (\Mg^{\rm ad}_{F_{(2)}})\to \pmb{\mu}_{2,F_{(2)}}\to 1$$ 
(to be compared with [6, Vol. III, Exp. XXIV, Thm. 1.3]). We note that if $n\neq 4$, then ${\it Aut}^\prime(\Mg^{\rm ad}_{F_{(2)}})=\Mg^{\rm ad}_{F_{(2)}}$ and if $n=4$, then $\Mg^{\rm ad}_{F_{(2)}}$ is the identity component of ${\it Aut}^\prime (\Mg^{\rm ad}_{F_{(2)}})$ and the quotient group scheme ${\bf E}:={\it Aut}^\prime(\Mg^{\rm ad}_{F_{(2)}})/\Mg^{\rm ad}_{F_{(2)}}$ is \'etale and a form of $({\bf Z}/3{\bf Z})_{F_{(2)}}$. We also have a second short exact sequence 
$$1\to \pmb{\mu}_{2,F_{(2)}}\to \Mg_{F_{(2)}}\to \Mg^{\rm ad}_{F_{(2)}}\to 1.$$ 
\indent
The non-trivial torsors of ${\pmb{\mu}}_{2,F}\tilde\to {({\bf Z}/2{\bf Z})}_F$ correspond to quadratic field extensions of $F$. We easily get that there exists a smallest totally real field extension $F_{\rm in}$ of $F$ of degree at most $2$ and such that there exists a class $\gamma^\prime\in H^1(F_{\rm in},{\it Aut}^\prime(\Mg^{\rm ad}_{F_{\rm in}}))$ that defines the form $\Md_{n,\tau,F_{\rm in}}$ of $\Mg_{F_{\rm in}}$ (more precisely, $F_{\rm in}$ is the extension of $F$ defined by the image of $\gamma$ in  $H^1(F,\pmb{\mu}_{2,F})$ via the first short exact sequence). The extensions of $\Mg_{F_{(2)}}$ and $\Md_{n,\tau,F_{(2)}}$ via ${\bf Z}_{(2)}$-monomorphisms $F_{(2)}\hookrightarrow W({\bf F})$ are isomorphic to ${\bf SO}^+_{2n,W({\bf F})}$. This implies that the field $F_{\rm in}$ is unramified above all primes of $F$ that divide $2$. 
\smallskip
If $n=4$, then from [13, Cor. 2 of p. 182] we easily get that the image of $\gamma^\prime$ in $H^1(F_{\rm in},{\bf E}_{F_{\rm in}})$ is the trivial class. Thus for all $n\ge 2$, there exists a class $\gamma_1\in H^1(F_{\rm in},\Mg^{\rm ad}_{F_{\rm in}})$ that defines the inner form $\Md_{n,\tau,F_{\rm in}}$ of $\Mg_{F_{\rm in}}$ (i.e., the images of $\gamma_1$ and $\gamma^\prime$ in $H^1(F_{\rm in},{\it Aut}^\prime(\Mg^{\rm ad}_{F_{\rm in}}))$ coincide). Let $\gamma_2\in H^2(F_{\rm in},\pmb{\mu}_{2,F_{\rm in}})$ be the co-boundary of $\gamma_1$ with respect to the second short exact sequence. Let $M(\gamma_2)$ be the central semisimple $F_{\rm in}$-algebra that defines $\gamma_2$; it is either $F_{\rm in}$ itself or a non-trivial form of $M_2(F_{\rm in})$. We know that the class $\gamma_2$ becomes trivial under each embedding of $F_{\rm in}$ into either ${\bf R}$ or $B({\bf F})$. Thus from [7, Lem. 5.5.3] we get that there exists a maximal torus of the reductive group scheme of invertible elements of $M(\gamma_2)$ that is defined by a totally real number field extension $F_{\rm oin}$ of $F_{\rm in}$ of degree at most $2$ and unramified above the primes of $F_{\rm in}$ that divide $2$. Let $F_1$ be the Galois extension of ${\bf Q}$ generated by $F_{\rm oin}$; it is totally real and unramified above $2$. The image of $\gamma_2$ in $H^2(F_1,\pmb{\mu}_{2,F_1})$ is the trivial class and thus we can speak about the class 
$$\gamma_0\in H^1(F_{1,(2)},\Md_{n,\tau,F_{1,(2)}})=H^1({\bf Z}_{(2)},{\rm Res}_{F_{1,(2)}/{\bf Z}_{(2)}} \Md_{n,\tau,F_{1,(2)}})$$ 
that defines the inner twist $\Mg_{F_{1,(2)}}$ of $\Md_{n,\tau,F_{1,(2)}}$. 
\smallskip
Let $q:=[F_1:{\bf Q}]$; we have $q\in r{\bf N}$ and our notations match with the ones of Subsubsection 2.2.5. Let $\tilde\gamma_0$ be the image of $\gamma_0$ in $H^1({\bf Z}_{(2)},{\bf Sp}(\tilde L_{(2)},\tilde\psi))$ via the ${\bf Z}_{(2)}$-monomorphisms of (1). We define $\tilde W:=\tilde L_{(2)}[{1\over 2}]$. We check that the image $\tilde\gamma_{0,{\bf Q}}$ of $\tilde\gamma_0$ in $H^1({\bf Q},{\bf Sp}(\tilde W,\tilde\psi))$ is the trivial class. Based on [7, Main Thm.], it suffices to show that the image of $\gamma_0$ in $H^1({\bf R},{\bf Sp}(\tilde W\otimes_{{\bf Q}} {\bf R},\tilde\psi))$ is trivial; but this is so as the image of $\gamma_0$ in $H^1({\bf R},({\rm Res}_{F_{1,(2)}/{\bf Z}_{(2)}} \Md_{n,\tau,F_{1,(2)}})_{\bf R})$ is trivial (as we have seen that $\Md_{n,\tau,F_{(2)}}$ and $\Mg_{F_{(2)}}$ become isomorphic under extensions via ${\bf Z}_{(2)}$-monomorphisms $F_{(2)}\hookrightarrow{\bf R}$). As ${\tilde\gamma}_{0,\bf Q}$ is trivial, we get that $\tilde\gamma_0$ itself is the trivial class. Thus by twisting the ${\bf Z}_{(2)}$-monomorphisms of (1) via $\gamma_0$, we get ${\bf Z}_{(2)}$-monomorphisms of the form 
$${\rm Res}_{F_{1,(2)}/{\bf Z}_{(2)}} \Mg_{F_{1,(2)}}\hookrightarrow G_{{\bf Z}_{(2)}}^{\prime,\rm der}\hookrightarrow {\bf Sp}(\tilde L_{(2)},\tilde\psi),$$ where $G_{{\bf Z}_{(2)}}^{\prime,\rm der}$ is an inner  form of $\Md_{nq}$. Using the natural ${\bf Z}_{(2)}$-monomorphism 
$$G^{\rm der}_{{\bf Z}_{(2)}}={\rm Res}_{F_{(2)}/{\bf Z}_{(2)}} \Mg_{F_{(2)}}\hookrightarrow {\rm Res}_{F_{1,(2)}/{\bf Z}_{(2)}} \Mg_{F_{1,(2)}},$$ 
we end up with a sequence of closed embedding ${\bf Z}_{(2)}$-monomorphisms 
$$G^{\rm der}_{{\bf Z}_{(2)}}\hookrightarrow G_{{\bf Z}_{(2)}}^{\prime,\rm der}\hookrightarrow {\bf Sp}(\tilde L_{(2)},\tilde\psi).\leqno (2)$$ 
\indent
As in [19, Subsect. 3.5] we argue that the normal subgroup $G^\flat:=G\cap {\bf Sp}(W,\psi)$ of $G$ is connected and thus reductive. Let $G^\flat_{{\bf Z}_{(2)}}$ be the schematic closure of $G^\flat$ in $G_{{\bf Z}_{(2)}}$. As in [19, Subsect. 3.5 and Subsubsect. 3.5.1] we argue that we have $G^{\flat}_{{\bf Z}_{(2)}}=G^{\rm der}_{\bf Z_{(2)}}$.
Thus $G_{{\bf Z}_{(2)}}$ is the flat, closed subgroup scheme of ${\bf GL}_{\tilde L_{(2)}}$ generated by $G^{\rm der}_{{\bf Z}_{(2)}}$ and by the center of ${\bf GL}_{\tilde L_{(2)}}$. Let $G^\prime_{{\bf Z}_{(2)}}$ be the flat, closed subgroup scheme of ${\bf GL}_{\tilde L_{(2)}}$ generated by $G_{{\bf Z}_{(2)}}^{\prime,\rm der}$ and by the center of ${\bf GL}_{\tilde L_{(2)}}$; it is a form of ${\bf GSO}^+_{2nq,{\bf Z}_{(2)}}$ and thus a reductive group scheme for which the notations match (i.e., whose derived group scheme is $G_{{\bf Z}_{(2)}}^{\prime,\rm der}$). As the group $\Md^{F_{1,(2)}/{\bf Z}_{(2)}}_{nq,\tau}$ of Subsubsection 2.2.5 splits over ${\bf Z}_2$ (cf. Lemma 2.2.6 (c)) and as the class $\gamma_0$ has a trivial image in $H^1({\bf Z}_2,{\rm Res}_{F_{1,(2)}\otimes_{{\bf Z}_{(2)}}{\bf Z}_2/{\bf Z}_2} \Md_{n,\tau,F_{1,(2)}\otimes_{{\bf Z}_{(2)}} {\bf Z}_2})$ (cf. Lang's theorem applied  to ${\rm Res}_{F_{1,(2)}\otimes_{{\bf Z}_{(2)}}{\bf Z}_2/{\bf Z}_2}  \Md_{n,\tau,F_{1,(2)}\otimes_{{\bf Z}_{(2)}} {\bf Z}_2}$ and the fact that the ring ${\bf Z}_2$ is henselian), the group scheme $G^{\prime,\rm der}_{{\bf Z}_2}$ is the trivial form of $(\Md_{nq,\tau}^{F_{1,(2)}/{\bf Z}_{(2)}})_{{\bf Z}_2}$ and therefore it is split. Thus the extension of $G^{\prime}_{{\bf Z}_{(2)}}$ to ${\bf Z}_2$ also splits and therefore it is isomorphic to ${\bf GSO}^+_{2nq,{\bf Z}_2}$. Thus the property 1.2 (i) holds.
\medskip\smallskip\noindent
{\bf 3.4. The new Shimura pair $(G^\prime,\Mx^\prime)$.} We define $G^\prime:=G_{{\bf Z}_{(2)}}^\prime\times_{{\bf Z}_{(2)}} {\bf Q}$. Let $\Mx^\prime$ be the $G^\prime({\bf R})$-conjugacy class of the composite of any element $h:{\bf S}\hookrightarrow G_{{\bf R}}$ of $\Mx$ with the ${\bf R}$-monomorphism $G_{{\bf R}}\hookrightarrow G_{{\bf R}}^\prime$. The Lie monomorphism $G_{{\bf R}}^{\rm der}({\bf R})\hookrightarrow G_{{\bf R}}^{\prime,\rm der}({\bf R})$ can be identified with the composite of a diagonal Lie monomorphism ${\bf SO}^*(2n)^r\hookrightarrow {\bf SO}^*(2n)^q$ with the Lie monomorphism $e_{n,q}: {\bf SO}^*(2n)^q\hookrightarrow {\bf SO}^*(2nq)$ (cf. the constructions of Subsection 3.3 and Lemma 2.2.6 (d)). As $({\rm Ad}\circ h)({\bf R})(i)$ is a Cartan involution of ${\rm Lie}(G^{\rm ad}_{{\bf R}})$ (cf. beginning of Section 1), the image $S_h$ through $h$ of the ${\bf SO}(2)={\bf SO}^*(2)$ Lie subgroup of ${\bf S}({\bf R})$, is the center of a maximal compact Lie subgroup of ${\bf SO}^*(2n)^r=G^{\rm der}({\bf R})$. But all maximal compact Lie subgroups of ${\bf SO}^*(2n)^r$ are ${\bf SO}^*(2n)^r$-conjugate (see [8, Ch. VI, Sect. 2]). By combining the last two sentences with Lemma 2.2.2 (c), we get that the Lie subgroup $S_h$ of ${\bf SO}^*(2nq)=G^{\prime,\rm der}({\bf R})$ is ${\bf SO}^*(2nq)$-conjugate to $C_{1,nq}={\rm Im}(s_{1,nq})$. Thus the centralizer of $S_h$ in ${\bf SO}^*(2nq)=G^{\prime,\rm der}({\bf R})$ is a maximal compact Lie subgroup of ${\bf SO}^*(2nq)$ that is isomorphic to ${\bf U}(nq)$, cf. Lemma 2.2.2 (c). This implies that the inner conjugation through $h({\bf R})(i)$ is a Cartan involution of ${\rm Lie}(G^{\prime,\rm der}_{{\bf R}})={\rm Lie}(G^{\prime,\rm ad}_{{\bf R}})$, cf. the classification of Cartan involutions in [8, Ch. X, Sect. 2]. 
\smallskip
The representation of $G^{\rm der}_{{\bf C}}$ on $\tilde W\otimes_{{\bf Q}} {\bf C}$ is a direct sum of standard representations of dimension $2n$ of the ${\bf SO}_{2n,{\bf C}}$ factors of $G^{\rm der}_{{\bf C}}$. Thus the Hodge ${\bf Q}$--structure on $\tilde W$ defined by $h$ has the same type as the Hodge ${\bf Q}$--structure on $W$ defined by $h$ and thus it is of type $\{(-1,0),(0,-1)\}$. As $G^{\prime,\rm der}({\bf R})$ is isomorphic to ${\bf SO}^*(2nq)$, $G^{\prime,{\rm ad}}$ is a simple ${\bf Q}$--group whose extension to ${\bf R}$ is non-compact. Based on the last two sentences and on the last sentence of the previous paragraph, we get that Deligne's axioms of the first paragraph of Section 1 hold for the pair $(G^\prime,\Mx^\prime)$. Thus $(G^\prime,\Mx^\prime)$ is a Shimura pair. 
\medskip\smallskip\noindent
{\bf 3.5. End of the proof of the Basic Theorem.} Let $\widetilde{\got A}$ be the free ${\bf Z}_{(2)}$-module of alternating forms on $\tilde L_{(2)}$ fixed by $G_{{\bf Z}_{(2)}}^{\prime,\rm der}$. There exist elements of $\widetilde{\got A}\otimes_{{\bf Z}_{(2)}} {\bf R}$ that define polarizations of the Hodge ${\bf Q}$--structure on $\tilde W$ defined by a fixed element $h\in\Mx$, cf. [4, Cor. 2.3.3]. Thus the real vector space $\widetilde{\got A}\otimes_{{\bf Z}_{(2)}} {\bf R}$ has a non-empty, open subset of such polarizations (cf. [4, Subsubsect. 1.1.18 (a)]). A standard application to $\widetilde{\got A}$ of the approximation theory  for independent  valuations, implies the existence of $\tilde\psi^\prime\in \widetilde{\got A}$ that is congruent modulo $2{\bf Z}_{(2)}$ to $\tilde\psi$ and that defines a polarization of the Hodge ${\bf Q}$--structure on $\tilde W$ defined by $h\in\Mx$. As $\tilde\psi^\prime$ is congruent modulo $2{\bf Z}_{(2)}$ to $\tilde\psi$, it is a perfect and alternating form on $\tilde L_{(2)}$. By replacing $\tilde\psi$ by  $\tilde\psi^\prime$, we can assume that $\tilde\psi$ defines a polarization of the Hodge ${\bf Q}$--structure on $\tilde W$ defined by $h\in\Mx$.
\smallskip
We get injective maps 
$\tilde f:(G,\Mx)\hookrightarrow ({\bf GSp}(\tilde W,\tilde\psi),\tilde\Ms)$ and $$\tilde f^\prime:(G^\prime,\Mx^\prime)\hookrightarrow ({\bf GSp}(\tilde W,\tilde\psi),\tilde\Ms)$$ 
of Shimura pairs. Let $\tilde L$ be a ${\bf Z}$-lattice of $\tilde W$ such that $\tilde\psi$ induces a perfect and alternating form on it and we have $\tilde L_{(2)}=\tilde L\otimes_{{\bf Z}} {\bf Z}_{(2)}$. Let $\tilde\Mb$, $\tilde\Mb^\prime$, and $v^\prime$ be as in the property 1.2 (ii). Let $\tilde\Mi$ be the involution of ${\rm End}(\tilde L_{(2)})$ defined by $\tilde\psi$. We check that the axioms 1.1 (i) to (v) hold for the quadruple $(\tilde f,\tilde L,v,\tilde\Mb)$. Obviously the axiom 1.1 (v) holds. We know that the axiom 1.1 (iv) holds for $(\tilde f,\tilde L,v,\tilde\Mb)$, cf. the last paragraph of Subsection 3.3. From Lemma 2.2.6 (ii) we get that the representation of $G^{\rm der}_{W({\bf F})}=G^{\flat}_{W({\bf F})}$ on $\tilde L_{(2)}\otimes_{{\bf Z}_{(2)}} W({\bf F})$ is isomorphic to the direct sum of a finite number of copies of the representation $\rho_{n,W({\bf F})}$. As $n\ge 2$, the fibres of $\rho_{n,W({\bf F})}$ are absolutely simple representations. From the last two sentences we get that $\tilde\Mb\otimes_{{\bf Z}_{(2)}} W({\bf F})$ is a product of matrix $W({\bf F})$-algebras. Thus the axiom 1.1 (i) holds for $(\tilde f,\tilde L,v,\tilde\Mb)$.
\smallskip
As $G^{\rm der}_{F_1}$ is a product of groups that are forms of ${\bf SO}^+_{2n,F_1}$ and that are permuted transitively by ${\rm Gal}(F_1/{\bf Q})$, the ${\bf Q}$--algebra $\tilde\Mb[{1\over 2}]$ is simple. Thus the axiom 1.1 (ii) holds for $(\tilde f,\tilde L,v)$. The fact that the axiom 1.1 (iii) holds for $(\tilde f,\tilde L,v)$ is a standard consequence of the fact that $G$ is generated by the center of ${\bf GL}_{\tilde W}$ and by $G^{\rm der}$ and of the description of the representation of $G_{W({\bf F})}^{\rm der}$ on $\tilde L_{(2)}\otimes_{{\bf Z}_{(2)}} W({\bf F})$. We conclude that the quadruple $(\tilde f,\tilde L,v,\tilde\Mb)$ is a hermitian orthogonal standard PEL situation in mixed characteristic $(0,2)$. Similarly we argue that $(\tilde f^\prime,\tilde L,v^\prime,\tilde\Mb^\prime)$ is a hermitian orthogonal standard Hodge situation in mixed characteristic $(0,2)$. Thus the property 1.2 (ii) holds. From the construction of (2) and Lemma 2.2.6 (ii) and (iii) we get:
\medskip
{\bf (i)} the natural $W({\bf F})$-monomorphism $G_{W({\bf F})}^{\rm der}\hookrightarrow G_{W({\bf F})}^{\prime,\rm der}$ is the composite of a diagonal $W({\bf F})$-monomorphism $\Md^r_{n,W({\bf F})}\hookrightarrow \Md^q_{n,W({\bf F})}$ with a standard $W({\bf F})$-monomorphism $d_{n,q,W({\bf F})}:\Md^q_{n,W({\bf F})}\hookrightarrow \Md_{nq,W({\bf F})}$;
\smallskip
{\bf (ii)} the faithful representation $G_{W({\bf F})}^{\prime,\rm der}\hookrightarrow {\bf GL}_{\tilde L_{(2)}\otimes_{\bf Z_{(2)}} W({\bf F})}$ is isomorphic to the direct sum of four copies of the representation $\rho_{nq,W({\bf F})}$. 
\medskip
From properties (i) and (ii) and from Lemma 2.1.1 we get that $G_{W({\bf F})}^{\rm der}$ is the identity component of the subgroup scheme of $G_{W({\bf F})}^{\prime,\rm der}$ that centralizes $\tilde\Mb\otimes_{\bf Z_{(2)}} W({\bf F})$. This implies that the property 1.2 (iii) also holds. This ends the proof of the Basic Theorem.\endproof  
\bigskip\smallskip
\centerline{\bigsll {\bf 4. Basic crystalline properties}}
\bigskip\smallskip
We will use the notations of the list $({\sharp}_1)$ of Section 3 and of the new list of notations 
$$n,\,r,\,d,\,\Mb_1,\,G_2,\,G_{2,{\bf Z}_{(2)}},\,B({\bf F}),\,F,\,\Mg,\,\kappa,\,j\in\kappa,\,F_j,\,F_{(2)},\,\Mg_{F_{(2)}},O_j\leqno ({\sharp}_2)$$
introduced in Section 3. Let $(\Ma,\Lambda)$
be the pull-back to $\Mn$ of the universal principally polarized abelian scheme over $\Mm$. Let $k$ be an arbitrary algebraically closed field of characteristic $2$ and of countable transcendental degree. We fix a ${\bf Z}_{(2)}$-embedding $O_{(v)}\hookrightarrow W(k)$ into the ring of $2$-typical Witt vectors with coefficients in $k$; all pull-backs to either ${\rm Spec}\, W(k)$ or ${\rm Spec}\, k$ of $O_{(v)}$-schemes, will be via it. Let $\sigma$ be the Frobenius automorphism of $W(k)$. For a $W(k)$-morphism $y:{\rm Spec}\, k\to\Mn_{W(k)}$ let
$$(A,\lambda_A):=y^*((\Ma,\Lambda)\times_{\Mn} \Mn_{W(k)}).$$ 
Let $(M,\phi,\psi_M)$ be the principally quasi-polarized (contravariant) Dieudonn\'e module over $k$ of the principally quasi-polarized $2$-divisible group of $(A,\lambda_A)$. Thus $M$ is a free $W(k)$-module of finite rank, $\phi:M\to M$ is a $\sigma$-linear endomorphism such that we have $2M\subset\phi(M)$, and $\psi_M$ is a perfect and alternating form on $M$ such that we have $\psi_M(\phi(x),\phi(u))=2\sigma(\psi_M(x,u))$ for all $x,u\in M$. We denote also by $\psi_M$ the perfect and alternating form on $M^\vee$ induced naturally by $\psi_M$. 
\smallskip
For $b\in\Mb$, we denote also by $b$ the ${\bf Z}_{(2)}$-endomorphism of $\Ma$ defined naturally by $b$ (cf. [19, Subsubsect. 4.1.1]). We denote also by $b$ different de Rham (crystalline) realizations of ${\bf Z}_{(2)}$-endomorphisms that correspond to $b$. Thus we will speak about the ${\bf Z}_{(2)}$-monomorphism $\Mb^{\rm opp}\hookrightarrow {\rm End}(M)$ that makes $M$ to be a $\Mb^{\rm opp}\otimes_{{\bf Z}_{(2)}} W(k)$-module and $M^\vee$ to be a $\Mb\otimes_{{\bf Z}_{(2)}} W(k)$-module; here $\Mb^{\rm opp}$ is the opposite ${\bf Z}_{(2)}$-algebra of $\Mb$. In this section we recall basic crystalline properties that are (or are proved) as in [19]. 
\medskip\smallskip\noindent
{\bf 4.1. Extra tensors.} Let $(v_\al)_{\al\in\Mj}$ be a family of tensors of $\Mt(W)$ such that $G$ is the subgroup of ${\bf GL}_W$ that fixes $v_{\al}$ for all $\al\in\Mj$, cf. [5, Prop. 3.1 (c)]. We choose the set $\Mj$ such that $\Mb\subseteq\Mj$ and for $b\in\Mb$ we have $v_b=b\in {\rm End}(W)=W\otimes_{{\bf Q}} W^\vee$.
\smallskip
Let $(\Mb_1\otimes_{{\bf Z}_{(2)}} {\bf Z}_2,\Mi)=\oplus_{j\in \kappa} (\Mb_j,\Mi)$ be the product decomposition of $(\Mb_1\otimes_{{\bf Z}_{(2)}} {\bf Z}_2,\Mi)$ into simple factors. Each $\Mb_j$ is a two sided ideal of $\Mb_1\otimes_{{\bf Z}_{(2)}} {\bf Z}_2$ that is a simple ${\bf Z}_2$-algebra and whose center is the ring of integers $O_j$ of the $2$-adic field $F_j$. As $\Mi$ is of first orthogonal type, it fixes the center of each $\Mb_j$.
\smallskip
As $F_j$ is unramified over ${\bf Q}_2$ (see Subsection 3.2), $O_j$ is a finite, \'etale ${\bf Z}_2$-algebra. We can identify $\Mb_j$ with ${\rm End}(\Mv_j)$, where $\Mv_j$ is a free $O_j$-module of rank $2n$. Let $s_j\in{\bf N}$ be such that as $\Mb_1\otimes_{{\bf Z}_{(2)}} {\bf Z}_2$-modules we can identify
$$L_{(2)}\otimes_{{\bf Z}_{(2)}} {\bf Z}_2=\oplus_{j\in\kappa} \Mv^{s_j}_j.\leqno (3)$$
We have $s_j\ge 2$, as the representation of $G^{\rm der}_{{\bf Z}_2}=G^{\flat}_{{\bf Z}_2}$ on $L_{(2)}\otimes_{{\bf Z}_{(2)}} {\bf Z}_2$ is symplectic. As $G^{\rm ad}$ is a simple ${\bf Q}$--group, the product $s_j[F_j:{\bf Q}_2]$ does not depend on $j\in\kappa$. We can redefine the direct summand  $\Mv^{s_j}_j$ of $L_{(2)}\otimes_{{\bf Z}_{(2)}} {\bf Z}_2$ as the maximal ${\bf Z}_2$-submodule that is generated by non-trivial, simple ${\rm Res}_{O_j/{\bf Z}_2} (\Mg_{F_{(2)}}\times_{F_{(2)}} O_j)$-submodules (we recall that ${\rm Res}_{O_j/{\bf Z}_2} (\Mg_{F_{(2)}}\times_{F_{(2)}} O_j)$ is a direct factor of $G^{\rm der}_{{\bf Z}_2}$ introduced in Subsection 3.2). 
\smallskip
Let $b_j$ be a perfect bilinear form on the $O_j$-module $\Mv_j$ that defines the involution $\Mi$ of $\Mb_j$, cf. [19, Lem. 3.3.1 (a)]. Thus $b_j$ is unique up to a ${{\bf G}_m}(O_j)$-multiple (cf. [19, Lem. 3.3.1 (b)]), it is fixed by $G^{\rm der}_{{\bf Z}_2}=G^{\flat}_{{\bf Z}_2}$, and it is symmetric (as $(\Mb_j,\Mi)$ is of orthogonal first type). Let $b_0:=\oplus_{j\in \kappa} b_j^{s_j}$; it is a perfect, symmetric bilinear form on the ${\bf Z}_2$-module $L_{(2)}\otimes_{{\bf Z}_{(2)}} {\bf Z}_2$ fixed by $G^{\rm der}_{{\bf Z}_2}=G^{\flat}_{{\bf Z}_2}$. The ${\bf Z}_2$-span of either $b_j$ or $b$ is normalized by $G_{{\bf Z}_2}$.
\medskip\smallskip\noindent
{\bf 4.2. Lifts of $y$.} As in [19, Subsect. 4.1] we argue that $\Mn_{E(G,\Mx)}$ is a closed subscheme of $\Mm_{E(G,\Mx)}$ and that there exists a compact, open subgroup $H_0$ of $G({\bf A}_f^{(2)})$ such that it acts freely on $\Mm$ and thus also on $\Mn$,  the schematic closure $\Mn/H_0$ of $\Mn_{E(G,\Mx)}/H_0={\rm Sh}(G,\Mx)/(H_2\times H_0)$ in $\Mm_{O_{(v)}}/H_0$ is a quasi-projective, normal $O_{(v)}$-scheme of finite type, and moreover $\Mn$ is a pro-\'etale cover of $\Mn/H_0$. The flat, finite type morphism $\Mn_{W(k)}/H_0\to {\rm Spec}\, W(k)$ has quasi-sections whose images contain the $k$-valued point of $\Mn_{W(k)}/H_0$ defined by $y$ (cf. [19, Subsubsect. 4.1.1]). This implies that the $W(k)$-morphism $y:{\rm Spec}\, k\to\Mn_{W(k)}$ has a lift 
$$z:{\rm Spec}\, V\to\Mn_{W(k)},$$
where $V$ is a finite, discrete valuation ring extension of $W(k)$. We define $(A_V,\lambda_{A_V}):=z^*((\Ma,\Lambda)_{\Mn_{W(k)}}).$ Let $B(k):=W(k)[{1\over 2}]$. As in [19, Subsect. 4.2] we argue that:
\medskip
{\bf (a)} for each $\al\in\Mj$ there exists a tensor $t_{\al}\in\Mt(M[{1\over 2}])$ that correspond naturally to $v_{\al}$ via Fontaine comparison theory for $A_V$;\smallskip
{\bf (b)} there exists a $B(k)$-linear isomorphism  $j_y:L_{(2)}\otimes_{{\bf Z}_{(2)}} B(k)\tilde\to M^\vee[{1\over 2}]$ 
that takes $\psi$ to $\psi_M$ and takes $v_{\al}$ to $t_{\al}$ for all $\al\in\Mj$. 
\medskip
Let $J$ be the schematic closure in ${\bf GL}_M$ of the subgroup $J_{B(k)}$ of ${\bf GL}_{M[{1\over 2}]}$ that fixes $t_{\al}$ for all $\al\in\Mj$. Equivalently, let $J$ be the schematic closure in ${\bf GSp}(M,\psi_M)$ of the identity component of the subgroup of ${\bf GSp}(M[{1\over 2}],\psi_M)$ that fixes $t_b$ for all $b\in\Mb\subseteq\Mj$. Thus, due to the property 1.1 (iii), the group $J_{B(k)}^{\rm der}$ is the identity component of the subgroup of ${\bf GL}_{M[{1\over 2}]}$ that fixes the subset $\Mb^{\rm opp}$ of ${\rm End}(M[{1\over 2}])$ and the involution of ${\rm End}(M[{1\over 2}])$ defined by $\psi_M$. The existence of $j_y$ implies that $J_{B(k)}$ is isomorphic to $G_{B(k)}$ and thus it is a reductive group. Let $j\in\kappa$. Each projection of $L_{(2)}\otimes_{{\bf Z}_{(2)}} {\bf Z}_2$ on a factor $\Mv_j$ of (3) along the direct sum of the other factors of (3), is an element of $\Mb\otimes_{{\bf Z}_{(2)}} {\bf Z}_2$. Thus to (3) corresponds naturally a direct sum decomposition
$$(M,\phi)=\oplus_{j\in\kappa} (N_j,\phi)^{s_j}\leqno (4) $$ 
of Dieudonn\'e modules over $k$. As the ${\bf Z}_2$-span of $b_j$ is normalized by $G_{\bf Z_2}$, we can speak about the perfect, symmetric bilinear form $c_j$ on $N_j$ that is the crystalline realization of $b_j$; it is uniquely determined up to a unit of ${\bf Z}_2$. As $\Mv_j$ is an $O_j$-module, $N_j$ is naturally a $W(k)\otimes_{{\bf Z}_2} O_j$-module and $c_j$ is $W(k)\otimes_{{\bf Z}_2} O_j$-linear. We consider the following condition:
\medskip
{\bf $(\natural)$} For each $j\in\kappa$, $c_j$ modulo $2W(k)$ is a perfect and alternating form on $N_j/2N_j$.
\medskip\noindent
{\bf 4.2.1. Proposition.} {\it We assume that the condition 4.2 $(\natural)$ holds. Then the flat, closed subgroup scheme $J$ of ${\bf GL}_M$ is reductive.}
\medskip\noindent
{\it Proof:} The formula $q_j(x):={{c_j(x,x)}\over 2}$ defines a quadratic form on the $W(k)\otimes_{{\bf Z}_2} O_j$-module $N_j$. The closed subgroup scheme ${\bf SO}(N_j,q_j)$ over $W(k)\otimes_{{\bf Z}_2} O_j=W(k)^{[F_j:\bf Q_2]}$ is isomorphic to ${\bf SO}^+_{2n,W(k)\otimes_{{\bf Z}_2} O_j}$ (cf. [19, Prop. 3.4 (c)]) and thus it is reductive. This implies that the schematic closure $J^{\rm der}$ of $J^{\rm der}_{B(k)}$ in ${\bf GL}_{M}$ is $\prod_{j\in\kappa} {\rm Res}_{W(k)\otimes_{{\bf Z}_2} O_j/W(k)} {\bf SO}(N_j,q_j)$ and thus it is a semisimple group scheme over $W(k)$ isomorphic to $({\bf SO}^+_{2n,W(k)})^r$ (we have $\sum_{j\in\kappa} [F_j:{\bf Q}_2]=[F:{\bf Q}]=r)$. If $Z$ is the center of ${\bf GL}_M$, then the intersection $Z\cap J^{\rm der}$ is a $\pmb{\mu}_{2,W(k)}$ group scheme. Let $J^\prime$ be the quotient of $Z\times_{W(k)} J^{\rm der}$ by a diagonal $\pmb{\mu}_{2,W(k)}$ closed subgroup scheme; it is a reductive group scheme (cf. [6, Vol. III, Exp. XXII, Prop. 4.3.1]) and we have a natural homomorphism $J^\prime\to {\bf GL}_M$ whose fibres are closed embeddings. The homomorphism $J^\prime\to {\bf GL}_M$ is a monomorphism (cf. [6, Vol. I, Exp. VI${}_B$, Cor. 2.11]) and thus also a closed embedding (cf. [6, Vol. II, Exp. XVI, Cor. 1.5 a)]) whose generic fibre can be identified with the closed embedding $J_{B(k)}\hookrightarrow {\bf GL}_{M[{1\over 2}]}$. Thus we can identify $J=J^\prime$ and therefore $J$ is a reductive group scheme over $W(k)$.\endproof
\medskip\noindent
{\bf 4.2.2. Proposition.} {\it We assume that the condition 4.2 $(\natural)$ holds. Then there exists a cocharacter $\mu:{\bf G}_{m,W(k)}\to J$ and a direct sum decomposition $M=F^1\oplus F^0$ such that the following two properties hold:
\medskip
{\bf (i)} each $\be\in {\bf G}_m(W(k))$ acts on $F^i$ as the multiplication with $\be^{-i}$ (here $i\in\{0,1\}$);
\smallskip
{\bf (ii)} the $k$-module $F^1/2F^1$ is the kernel of the reduction $\vph$ modulo $2W(k)$ of $\phi$.
\medskip
Moreover, the normalizer of $F^1/2F^1$ in the special fibre $J_k$ of $J$ is a parabolic subgroup $P_k$ of $J_k$ and the dimension of $\dim(J_k/P_k)$ is $d={{rn(n-1)}\over 2}$.}
\medskip\noindent
{\it Proof:} We have a direct sum decomposition ${\rm Ker}(\vph)=\oplus_{j\in\kappa} ((N_j/2N_j)\cap {\rm Ker}(\vph))^{s_j}$ and each intersection $(N_j/2N_j)\cap {\rm Ker}(\vph)$ is naturally a $k\otimes_{{\bf F}_2} O_j/2O_j$-module. Thus there exists a direct summand $\tilde F^1_j$ of $N_j$ that is a $W(k)\otimes_{{\bf Z}_2} O_j$-submodule and that lifts $(N_j/2N_j)\cap {\rm Ker}(\vph)$. Based on Proposition 4.2.1, the rest of the proof of this Proposition is the same as of [19, Prop. 6.1 and Cor. 6.1.1]. In other words, as in loc. cit. one first checks that: 
\medskip
{\bf (iii)} our hypothesis on $c_j$ implies that we have $c_j(\tilde F^1_j,\tilde F^1_j)\in 4W(k)\otimes_{{\bf Z}_2} O_j$ and
\smallskip
{\bf (iv)} property (iii) and [19, Prop. 3.4 (b)] imply that there exists a direct sum decomposition $N_j=F^1_j\oplus F^0_j$ of $W(k)\otimes_{{\bf Z}_2} O_j$-modules such that $F^1_j/2F^1_j=\tilde F^1_j/2\tilde F^1_j$ and we have $c_j(F^1_j,F^1_j)=c_j(F^0_j,F^0_j)=0$.
\medskip
We take $F^1:=\oplus_{j\in\kappa} (F^1_j)^{s_j}$ and $F^0:=\oplus_{j\in\kappa} (F^0_j)^{s_j}$ and $\mu$ to be defined by the property (i). As $\mu$ acts as scalar multiplication on $\oplus_{j\in\kappa} c_j^{s_j}$ and as we can identify $J^{\rm der}=\prod_{j\in\kappa} {\rm Res}_{W(k)\otimes_{{\bf Z}_2} O_j/W(k)} {\bf SO}(N_j,q_j)$ (cf. proof of Proposition 4.2.1 where $q_j$ is defined), we easily get that $\mu$ factors through $J$.  Obviously (ii) holds.\endproof
\medskip\noindent
{\bf 4.2.3. Proposition.} {\it We assume that the condition 4.2 $(\natural)$ holds. Then there exist isomorphisms $L_{(2)}\otimes_{{\bf Z}_{(2)}} W(k)\tilde\to M^\vee$ of $\Mb\otimes_{{\bf Z}_{(2)}} W(k)$-modules that induce symplectic isomorphisms $(L_{(2)}\otimes_{{\bf Z}_{(2)}} W(k),\psi)\tilde\to (M^\vee,\psi_M)$ and that take $v_{\alpha}$ to $t_{\alpha}$ for all $\alpha\in\Mj$.}
\medskip\noindent
{\it Proof:} We refer to the $B(k)$-linear isomorphism  $j_y:L_{(2)}\otimes_{{\bf Z}_{(2)}} B(k)\tilde\to M^\vee[{1\over 2}]$ of the property 4.2 (b). Let $L_y:=j_y^{-1}(M^\vee)$. It is a $W(k)$-lattice of $L_{(2)}\otimes_{{\bf Z}_{(2)}} W(k)$ such that the following three properties hold:
\medskip
{\bf (i)} for all $b\in\Mb\otimes_{{\bf Z}_{(2)}} W(k)$ we have $b(L_y)\subseteq L_y$;
\smallskip
{\bf (ii)} the schematic closure of $G_{B(k)}$ in ${\bf GL}_{L_y}$ is a reductive group scheme $j_y^{-1} J j_y$ over $W(k)$ (cf. Proposition 4.2.1);
\smallskip
{\bf (iii)} we get a perfect and alternating form $\psi:L_y\otimes_{W(k)} L_y\to W(k)$ whose $W(k)$-span is normalized by $j_y^{-1} J j_y$. 
\medskip
As in [19, Subsect. 5.2] we argue that the properties (i) to (iii) imply that there exists an element $g\in G^{\flat}(B(k))$ such that we have $g(L_{(2)}\otimes_{{\bf Z}_{(2)}} W(k))=L_y$. By replacing $j_y$ with $j_yg$, we can assume that $j_y(L_{(2)}\otimes_{{\bf Z}_{(2)}} W(k))=j_y(L_y)=M^\vee$. Thus $j_y:L_{(2)}\otimes_{{\bf Z}_{(2)}} W(k)\tilde\to M^\vee$ is an isomorphism of $\Mb\otimes_{{\bf Z}_{(2)}} W(k)$-modules that induces a symplectic isomorphism $(L_{(2)}\otimes_{{\bf Z}_{(2)}} W(k),\psi)\tilde\to (M^\vee,\psi_M)$ and that takes $v_{\alpha}$ to $t_{\alpha}$ for all $\alpha\in\Mj$.\endproof 
\bigskip\smallskip
\centerline{\bigsll {\bf 5. Proof of the Main Theorem}}
\bigskip\smallskip
We will combine Sections 3 and 4 to prove the Main Theorem. We will use the notations of the lists $({\sharp}_1)$ and $({\sharp}_2)$ of Sections 3 and 4 (respectively). Also the notations 
$$\tilde f:(G,\Mx)\to ({\bf GSp}(\tilde W,\tilde\psi),\tilde\Ms),\,\tilde f^\prime:(G^\prime,\Mx^\prime)\to ({\bf GSp}(\tilde W,\tilde\psi),\tilde\Ms),\,\tilde L_{(2)},\,\tilde L,\,\tilde\Mb,\,\tilde\Mb^\prime,\,v^\prime\leqno ({\sharp}_3)$$will be as in Subsections 3.3 to 3.5. Let the field $k$ be as in Section 4. Let $\break\tilde K_2:={\bf GSp}(\tilde L_{(2)},\tilde\psi)({\bf Z}_2)$ and $H_2^\prime:=G^\prime({\bf Q}_2)\cap\tilde K_2=G_{{\bf Z}_{(2)}}^{\prime}({\bf Z}_2)$. Let $\tilde\Mn$, $(\tilde\Ma,\tilde\Lambda)$ (resp. $\tilde\Mn^\prime$, $(\tilde\Ma^\prime,\tilde\Lambda^\prime)$), and $\tilde\Mm$ be the analogues of $\Mn$, $(\Ma,\Lambda)$, and $\Mm$ but obtained working with the hermitian orthogonal standard PEL situation $(\tilde f,\tilde L,v,\tilde\Mb)$ (resp. $(\tilde f^\prime,\tilde L,v^\prime,\tilde\Mb^\prime)$) instead of with $(f,L,v,\Mb)$. Note that $\tilde\Mn$ is a closed subscheme of $\tilde\Mn^{\prime}_{O_{(v)}}$. We know that $\tilde\Mn^{\prime}$ is regular and formally smooth over $O_{(v^{\prime})}$, cf. [19, Thm. 1.4 (b)]. 
\smallskip
It is well known that the $E(G,\Mx)$-scheme $\tilde\Mn_{E(G,\Mx)}=\Mn_{E(G,\Mx)}$ is regular and formally smooth. We choose $H_0$ such that $\tilde\Mn$ is as well a pro-\'etale cover of the quasi-projective, normal $O_{(v)}$-scheme $\tilde\Mn/H_0$ of relative dimension $d={{rn(n-1)}\over 2}$.  Here $\tilde\Mn/H_0$ is defined as well to be the schematic closure of $\tilde\Mn_{E(G,\Mx)}$ in $\tilde\Mm_{O_{(v)}}/H_0$. 
\medskip\noindent
{\bf 5.1. Theorem.} {\it The $O_{(v)}$-scheme $\tilde\Mn$ is regular and formally smooth.}
\medskip\noindent
{\it Proof:} As the $E(G,\Mx)$-scheme $\tilde\Mn_{E(G,\Mx)}=\Mn_{E(G,\Mx)}$ is regular and formally smooth, to prove the theorem it suffices to show that the $O_{(v)}$-scheme $\tilde\Mn/H_0$ is smooth at all points of positive characteristic. Thus it suffices to show that for each point $\tilde y\in\tilde\Mn_{W(k)}(k)$, the tangent space ${\bf T}_{\tilde y}$ of $\tilde y$ in $\tilde\Mn_k$ is a $k$-vector space of dimension $d$. 
\smallskip
We denote also by $\tilde y$ the $k$-valued point of $\tilde\Mn_{W(k)}$ defined by $\tilde y$. Let $(\tilde M,\tilde\phi,\psi_{\tilde M})$ be the principally quasi-polarized (contravariant) Dieudonn\'e module over $k$ of the principally quasi-polarized $2$-divisible group of 
$$(\tilde A,\lambda_{\tilde A}):=\tilde y^*((\tilde\Ma,\tilde\Lambda)\times_{\tilde\Mn} \tilde\Mn_{W(k)}).$$ 
Let 
$$\tilde L_{(2)}\otimes_{{\bf Z}_{(2)}} {\bf Z}_2=\oplus_{j\in\kappa} \tilde\Mv_j^{\tilde s_j}\leqno (5)$$ 
be the direct sum decomposition that is analogous to (3) (here each $\tilde s_j\in {\bf N}\setminus\{1\}$). This makes sense as (cf. end of Subsection 3.2) we can identify $\kappa$ with the set of factors of $G^{\rm der}_{{\bf Z}_2}$ that are Weil restrictions of semisimple ${\bf SO}_{2n}$ group schemes. Thus we can assume that the direct factor ${\rm Res}_{O_j/{\bf Z}_2} (\Mg_{F_{(2)}}\times_{F_{(2)}} O_j)$ of $G^{\rm der}_{{\bf Z}_2}$ acts non-trivially on $\tilde\Mv_j$. Moreover, $\tilde\Mv_j$ is a free $O_j$-module of rank $2n$ and is a $\Mg_{F_{(2)}}\times_{F_{(2)}} O_j$-module whose fibres are simple modules. The representations $\Mg_{F_{(2)}}\times_{F_{(2)}} O_j\to {\bf GL}_{\Mv_j}$ and $\Mg_{F_{(2)}}\times_{F_{(2)}} O_j\to {\bf GL}_{\tilde\Mv_j}$ over $O_j$ are isomorphic. This is so as over $W({\bf F})$ they are isomorphic to $\rho_{n,W({\bf F})}$ and as their fibres are absolutely irreducible (as $n\ge 2$). We conclude that:
\medskip
$\bullet$ The $G^{\rm der}_{{\bf Z}_2}$-modules (and thus also the $G_{{\bf Z}_2}$-modules) $\Mv_j$ and $\tilde\Mv_j$ are isomorphic.
\medskip
Each perfect, symmetric bilinear form fixed by $G^{\prime,{\rm der}}_{{\bf Z}_2}$ and obtained in the same way $b_0$ was (but working with the $G^{\prime,{\rm der}}_{{\bf Z}_2}$-module $\tilde L_{(2)}\otimes_{{\bf Z}_{(2)}} {\bf Z}_2$ instead of the $G^{{\rm der}}_{{\bf Z}_2}$-module $L_{(2)}\otimes_{{\bf Z}_{(2)}} {\bf Z}_2$), has the property that the isomorphism $\oplus_{j\in\kappa} \tilde\Mv_j^{\tilde s_j}\to \oplus_{j\in\kappa} (\tilde\Mv_j^{\tilde s_j})^{\vee}$ induced naturally by it maps the direct summand $\tilde\Mv_j^{\tilde s_j}$ onto the direct summand $(\tilde\Mv_j^{\tilde s_j})^{\vee}$ for all $j\in\kappa$. Based on this, it is easy to see that we can choose the direct sum decomposition (5) so that for each $j\in \kappa$, the isomorphism $\tilde\Mv_j^{\tilde s_j}\to (\tilde\Mv_j^{\tilde s_j})^{\vee}$ of $G^{\rm der}_{\bf Z_2}$-modules we have just obtained is the one induced by $\tilde b_j^{\tilde s_j}$, where $\tilde b_j$ is a perfect, symmetric bilinear form on $\tilde\Mv_j$ that is fixed by $G_{\bf Z_2}^{\rm der}$ and that is the natural analogue of $b_j$. In other words, we can choose a perfect, symmetric bilinear form $\tilde b_0:=\oplus_{j\in\kappa} \tilde b_j^{\tilde s_j}$ on $\tilde L_{(2)}\otimes_{{\bf Z}_{(2)}} {\bf Z}_2$ that is the analogue of $b_0=\oplus_{j\in\kappa} b_j^{s_j}$ of Subsection 4.1 and that is fixed by $G^{\prime,{\rm der}}_{{\bf Z}_2}$ as well. 
\smallskip
Let $(\tilde M,\tilde\phi)=\oplus_{j\in\kappa} (\tilde N_j^{\tilde s_j},\tilde\phi)$ be the direct sum decomposition that is the analogue of (4). Let $\tilde c_j$ be the perfect, symmetric bilinear form on $\tilde N_j$ that is the analogue of $c_j$ of Subsection 4.1. Let $\tilde J$ (resp. $\tilde J^\prime$) be the flat, closed subgroup scheme of ${\bf GL}_{\tilde M}$ that is the analogue of $J$ of Subsection 4.2 and that corresponds to $\tilde y\in\tilde\Mn_{W(k)}(k)$ (resp. to the point $\tilde y^\prime\in\tilde\Mn^{\prime}_{W(k)}(k)$ defined by $\tilde y$).  We know that $\tilde J^\prime$ is a reductive group scheme isomorphic to ${\bf GSO}^+_{2nq,W(k)}$, cf. [19, Subsect. 5.2] applied $\tilde y^\prime\in\tilde\Mn^{\prime}_{W(k)}(k)$ (based on the fact that Theorem 1.2 (i) holds). Let $\tilde c_0=\oplus_{j\in\kappa} \tilde c_j^{\tilde s_j}$ be the perfect, symmetric form on $\tilde M$ that corresponds naturally to $\tilde b_0$. 
\smallskip
As $G^{\prime,{\rm der}}_{{\bf Z}_2}$ fixes $\tilde b_0$, from [19, Prop. 5.1] applied to the point $\tilde y^\prime\in\tilde\Mn^{\prime}_{W(k)}(k)$, we get that $\tilde c_0$ modulo $2W(k)$ is alternating. Thus each $\tilde c_j$ modulo $2W(k)$ is alternating. Therefore $\tilde J$ is a reductive, closed subgroup scheme of ${\bf GL}_{\tilde M}$, cf. Proposition 4.2.1 applied to $\tilde y\in\tilde\Mn(k)$ instead of $y\in\Mn(k)$. Let $\tilde\mu:{\bf G}_{m,W(k)}\to\tilde J$ and $\tilde M=\tilde F^1\oplus\tilde F^0$ be the analogues of the cocharacter $\mu:{\bf G}_{m,W(k)}\to J$ and of the direct sum decomposition $M=F^1\oplus F^0$ introduced in Proposition 4.2.2. 
\smallskip
We have a natural direct sum decomposition into $W(k)$-modules
$${\rm End}(\tilde M)={\rm End}(\tilde F^1)\oplus {\rm End}(\tilde F^0)\oplus {\rm Hom}(\tilde F^1,\tilde F^0)\oplus {\rm Hom}(\tilde F^0,\tilde F^1)$$ as well as a modulo $2W(k)$ version of it. Let $\Mo_{\tilde y^\prime}$ and $\Mo^{\rm big}_{\tilde y^\prime}$ be the local rings of $\tilde y^\prime$ in $\tilde\Mn^{\prime}_{W(k)}$ and in $\tilde\Mm_{W(k)}$ (respectively). As we have a natural $W(k)$-epimorphism $\Mo^{\rm big}_{\tilde y^\prime}\twoheadrightarrow \Mo_{\tilde y^\prime}$, the tangent space ${\bf T}_{\tilde y^\prime}$ of  $\tilde y^\prime$ in $\tilde\Mn^{\prime}_k$ is naturally identified with the tensorization with $k$ of the image of the Kodaira--Spencer map of the natural pull-back of $\tilde\Ma$ to ${\rm Spec}\, \Mo_{\tilde y^\prime}$. In other words, ${\bf T}_{\tilde y^\prime}$ is naturally identified with the intersection ${\rm Lie}(\tilde J^\prime_k)\cap {\rm Hom}(\tilde F^1/2\tilde F^1,\tilde F^0/2\tilde F^0)$ taken inside ${\rm End}(\tilde M/2\tilde M)$ (cf. [19, Subsect. 6.4 and proof of Prop. 6.7]). Thus the tangent space ${\bf T}_{\tilde y}$ is a subspace of the intersection of ${\rm Lie}(\tilde J^\prime_k)\cap {\rm Hom}(\tilde F^1/2\tilde F^1,\tilde F^0/2\tilde F^0)$ with the centralizer of $\tilde\Mb^{\rm opp}\otimes_{\bf Z_{(2)}} k$ in ${\rm End}(\tilde M/2\tilde M)$. By applying Proposition 4.2.3 to $\tilde y$ and based on the property 1.2 (iii), we get that the identity component of the centralizer of $\tilde\Mb^{\rm opp}\otimes_{\bf Z_{(2)}} k$ in $\tilde J^\prime_k$ is $\tilde J_k$. Thus ${\bf T}_{\tilde y}$ is a $k$-vector subspace of the intersection ${\rm Lie}(\tilde J_k)\cap {\rm Hom}(\tilde F^1/2\tilde F^1,\tilde F^0/2\tilde F^0)$. 
\smallskip
Let $\tilde P_k$ be the parabolic subgroup of $\tilde J_k$ that is the normalizer of $\tilde F^1/2\tilde F^1$ in $\tilde J_k$, cf. Proposition 4.2.2 applied to $\tilde y$. As $\tilde\mu$ is a cocharacter of $\tilde J$, we have a direct sum decomposition into $k$-vector spaces ${\rm Lie}(\tilde J_k)={\rm Lie}(\tilde P_k)\oplus ({\rm Lie}(\tilde J_k)\cap {\rm Hom}(\tilde F^1/2\tilde F^1,\tilde F^0/2\tilde F^0))$.  Thus $\dim_k({\bf T}_{\tilde y})\le \dim_k({\rm Lie}(\tilde J_k)/{\rm Lie}(\tilde P_k))=\dim(\tilde J_k/\tilde P_k)$ and therefore (cf. Proposition 4.2.2 applied to $\tilde y$) we have $\dim_k({\bf T}_{\tilde y})\le d$. As we obviously have $\dim_k({\bf T}_{\tilde y})\ge d$, we conclude that $\dim_k({\bf T}_{\tilde y})=d$.\endproof
\medskip\smallskip\noindent
{\bf 5.2. Proposition.} {\it The natural identification of $\tilde\Mn_{E(G,\Mx)}$ with $\Mn^{\rm n}_{E(G,\Mx)}$ extends uniquely to an $O_{(v)}$-morphism $\Xi:\tilde\Mn\to \Mn^{\rm n}$.}
\medskip\noindent
{\it Proof:} To ease notations, let $\tilde Y:=\tilde\Mn_{W(k)}$; it is a regular scheme over $W(k)$ (cf. Theorem 5.1). Let $({\got D},\lambda_{\got D})$ and $(\tilde{\got D},\lambda_{\tilde{\got D}})$ be the principally quasi-polarized $2$-divisible groups of $(\Ma,\Lambda)_{\Mn^{\rm n}_{W(k)}}$ and  $(\tilde\Ma,\tilde\Lambda)_{\tilde Y}$ (respectively). To the decompositions (3) and (5) correspond naturally decompositions
$${\got D}=\oplus_{j\in\kappa} {\got D}_j^{s_j}\;\; {\rm and}\;\; \tilde{\got D}=\oplus_{j\in\kappa} \tilde{\got D}_j^{\tilde s_j}$$
(respectively) into $2$-divisible groups. The fact that the $G_{{\bf Z}_2}$-modules $\Mv_j$ and $\tilde\Mv_j$ are isomorphic (cf. proof of Theorem 5.1) can be encoded in the existence of a suitable ${\bf Z}_2$-endomorphism between the abelian schemes $\Ma_{\Mn_{B(k)}}$ and $\tilde\Ma_{\tilde\Mn_{B(k)}}$ over $\Mn_{B(k)}=\tilde Y_{B(k)}$; this  ${\bf Z}_2$-endomorphism allows us to identify naturally ${\got D}_{j,B(k)}$ with $\tilde{\got D}_{j,B(k)}$ as $2$-divisible groups over $\Mn^{\rm n}_{B(k)}=\tilde Y_{B(k)}$. Thus we can speak about the $2$-divisible group ${\got E}:=\oplus_{j\in\kappa} \tilde{\got D}_j^{s_j}$ over $\tilde Y$ which extends $\oplus_{j\in\kappa} \tilde{\got D}_{j,B(k)}^{s_j}={\got D}_{B(k)}$; note that both ${\got E}$ and ${\got D}$ are over $\tilde Y$  but ${\got E}$ involves the $s_j$' numbers while $\tilde{\got D}$ involves the $\tilde s_j$'s numbers. Let $\lambda_{{\got E},B(k)}:=\lambda_{{\got D},B(k)}$; it is a principal quasi-polarization of ${\got E}_{B(k)}={\got D}_{B(k)}$. As the $O_{(v)}$-scheme $\tilde Y$ is flat and normal, a theorem of Tate (see [15, Thm. 4]) implies that $\lambda_{{\got E},B(k)}$ extends uniquely to a principal quasi-polarization $\lambda_{\got E}$ of $\got E$.
\smallskip
Let $\tilde O$ be a local ring of $\tilde Y$ that is a discrete valuation ring of mixed characteristic $(0,2)$. Let $\tilde K$ be the field of fractions of $\tilde O$. We also view the natural $E(G,\Mx)$-morphism $\tilde z:{\rm Spec}\, \tilde K\to \tilde Y_{B(k)}$ as a $\tilde K$-valued point $\tilde z$ of $\Mn^{\rm n}_{B(k)}$. Thus we can speak about the abelian variety $\tilde A_{\tilde K}:=\tilde z^*(\Ma_{\Mn^{\rm n}_{B(k)}})$ over $\tilde K$. As $\tilde A_{\tilde K}$ has a level-$l$ structure for each odd natural number $l$, it extends to an abelian scheme $\tilde A_{\tilde O}$ over $\tilde O$ (cf. the N\'eron--Shafarevich--Ogg criterion of good reduction; see [1, Ch. 7, Sect. 7.4, Thm. 5]). It is easy to see that each principal quasi-polarization of  $\tilde A_{\tilde K}$ extends to a principal polarization of $\tilde A_{\tilde O}$ (at the level of isomorphisms between abelian scheme over $\tilde O$ this follows from [1, Ch. 7, Sect. 7.5, Prop. 3]). We conclude that there exists an open subscheme $\tilde U$ of the regular scheme $\tilde Y$ such that:
\medskip
{\bf (i)} we have $\tilde Y_{B(k)}\subseteq \tilde U$ and the codimension of $\tilde U$ in $\tilde Y$ is at least $2$;
\smallskip
{\bf (ii)} the principally polarized abelian scheme $(\Ma,\Lambda)_{\tilde Y_{B(k)}}=(\Ma,\Lambda)_{\Mn^{\rm n}_{B(k)}}$ extends to a principally polarized abelian scheme $(A_{\tilde U},\lambda_{A_{\tilde U}})$ over $\tilde U$.
\medskip
Tate theorem also implies that the principally quasi-polarized $2$-divisible group of $(A_{\tilde U},\lambda_{A_{\tilde U}})$ is $(\got E,\lambda_{\got E})_{\tilde U}$. As the scheme $\tilde Y$ is regular, from [18, proof of Prop. 4.1] we get that the principally polarized abelian scheme $(A_{\tilde U},\lambda_{A_{\tilde U}})$ over $\tilde U$ extends to a principally polarized abelian scheme $(A_{\tilde Y},\lambda_{A_{\tilde Y}})$ over $\tilde Y$ whose principally quasi-polarized $2$-divisible group is $(\got E,\lambda_{\got E})$. The natural level-$l$ symplectic similitude structures of $(\Ma,\Lambda)_{\tilde Y_{B(k)}}$ extend naturally to level-$l$ symplectic similitude structures of  $(A_{\tilde Y},\lambda_{A_{\tilde Y}})$. Thus the natural $E(G,\Mx)$-morphism $\tilde Y_{B(k)}\to \Mm_{B(k)}$ extends uniquely to a morphism $\tilde Y\to\Mm$. This last morphism factors through a $W(k)$-morphism $\Xi_{W(k)}:\tilde Y\to\Mn^{\rm n}_{W(k)}$ which is the pull-back to $W(k)$ of the searched for $O_{(v)}$-morphism $\Xi:\tilde\Mn\to \Mn^{\rm n}$ and for which we have $\Xi_{W(k)}^*({\got D})={\got E}$, cf. constructions. Obviously $\Xi$ is unique.\endproof
\medskip\noindent
{\bf 5.2.1. Remark.} The $O_{(v)}$-scheme $\Mn^{\rm n}$ has the extension property with respect to healthy regular $O_{(v)}$-schemes, cf. [17, Defs. 3.2.1 2) and 3.2.3 3), Ex. 3.2.9 and Cor. 3.4.1]. Thus the existence of $\Xi$ also follows from [20, Cor. 5] and Theorem 5.1 which imply that the $O_{(v)}$-scheme $\tilde\Mn$ is healthy regular.
\medskip\smallskip\noindent
{\bf 5.3. Lemma.} {\it We consider the $O_{(v)}$-morphism $\Xi_{H_0}:\tilde\Mn/H_0\to \Mn^{\rm n}/H_0$ whose generic fibre is a natural identification (isomorphism) and whose natural pull-back is the $O_{(v)}$-morphism $\Xi:\tilde\Mn\to\Mn^{\rm n}$. Then $\Xi_{H_0}$ is projective.}
\medskip\noindent
{\it Proof:} As $\Xi_{H_0}$ is quasi-projective and its generic fibre is an isomorphism, it suffices to show that for each discrete valuation ring $O$ that is a faithfully flat $O_{(v)}$-algebra, every $O_{(v)}$-morphism ${\rm Spec}\, O\to \Mn^{\rm n}/H_0$ factors uniquely through $\tilde\Mn/H_0$. To check this, we can assume that $O$ is strictly henselian and it suffices to show that each $O_{(v)}$-morphism $h_O:{\rm Spec}\, O\to \Mn^{\rm n}$ factors uniquely through $\tilde\Mn$. Let $K:=O[{1\over 2}]$. As in the proof of Theorem 5.1 we argue that the principally polarized abelian variety $(\tilde\Ma,\tilde\Lambda)_K$ over $K$ extends to a principally polarized abelian scheme over $O$ and that the generic fibre of $h_O$ extends uniquely to a morphism ${\rm Spec}\, O\to\tilde\Mm$ which factors through $\tilde\Mn$. Thus $h_O:{\rm Spec}\, O\to \Mn^{\rm n}$ factors uniquely through $\tilde\Mn$.\endproof
\medskip\noindent
{\bf 5.4. End of the proof of the Main Theorem.} As $\Xi_{W(k)}^*({\got D})={\got E}$ and as $\tilde{\got D}$ and thus also ${\got E}$ is a versal $2$-divisible group (i.e., its Kodaira--Spencer map is injective at all $k$-valued points of $\tilde Y$), the morphism $\Xi_{W(k)}$ and thus also $\Xi_{H_0}$ has finite fibres. Thus $\Xi_{H_0}$ is a quasi-finite morphism. From this and Lemma 5.3 we get that $\Xi_{H_0}$ is a finite, birational $O_{(v)}$-morphism. As the quasi-projective $O_{(v)}$-scheme $\Mn^{\rm n}/H_0$ is normal, we conclude that $\Xi_{H_0}$ is an isomorphism. Thus the $O_{(v)}$-morphism $\Xi:\tilde\Mn\to \Mn^{\rm n}$ is an isomorphism. Thus the $O_{(v)}$-scheme $\Mn^{\rm n}$ is regular and formally smooth (cf. Theorem 5.1) i.e., the Main Theorem holds.\endproof
\medskip\noindent
{\bf 5.5. Corollary.} {\it Let $y\in\Mn_{W(k)}(k)$. Let $(M,\phi,\psi_M)$ be as in the beginning of Section 4. Let $(v_{\alpha})_{\alpha\in\Mj}$ and $(t_{\alpha})_{\alpha\in\Mj}$ be as in Subsections 4.1 and 4.2 (respectively). Then there exist isomorphisms $L_{(2)}\otimes_{{\bf Z}_{(2)}} W(k)\tilde\to M^\vee$ of $\Mb\otimes_{{\bf Z}_{(2)}} W(k)$-modules that induce symplectic isomorphisms $(L_{(2)}\otimes_{{\bf Z}_{(2)}} W(k),\psi)\tilde\to (M^\vee,\psi_M)$ and that take $v_{\alpha}$ to $t_{\alpha}$ for all $\alpha\in\Mj$.}
\medskip\noindent
{\it Proof:} As the $O_{(v)}$-scheme $\Mn^{\rm n}$ is regular and formally smooth (cf. Main Theorem), there exists a lift $z\in \Mn_{W(k)}(W(k))$ of $y$. Let $F$ be the Hodge filtration of $M$ defined by $z^*(\Ma_{\Mn_{W(k)}})$. The decomposition (4) extends to a decomposition
$$(M,F,\phi)=\oplus_{j\in\kappa} (N_j,F_j,\phi)^{s_j}$$ 
of filtered $F$-crystals over $k$. For $x$, $u\in N_j$ we have $c_j(\phi(x),\phi(u))=2\sg(c_j(x,u))$. If $v\in F_j$, then $c_j(v,v)=0$ and $\phi(v)\in 2N_j$; thus $c_j({{\phi(v)}\over 2},{{\phi(v)}\over 2})=0$. As $N_j=\phi(N_j)+{1\over 2}\phi(F_j)$, each $x\in N_j$ is a sum ${1\over 2}\phi(v)+\phi(u)$, with $u\in N_j$ and $v\in F_j$. Thus 
$$c_j(x,x)=2c_j(\phi(u),{1\over 2}\phi(v))+c_j(\phi(u),\phi(u))=2c_j(\phi(u),{1\over 2}\phi(v))+2\sg(c_j(u,u))\in 2W(k)$$ 
i.e., $c_j$ modulo $2W(k)$ is alternating. Thus the corollary follows from Proposition 4.2.3.\endproof 
\medskip\noindent
{\it Acknowledgments.} We would like to thank University of Arizona, IAS--Princeton, and Binghamton University for good working conditions. We would also like to thank the referee for pointing out several typos. This research was partially supported by the NSF grant DMS \#0900967.
\bigskip\smallskip
\centerline{\bigsll {\bf References}}
\bigskip
\item{[1]} S. Bosch, W. L\"utkebohmert, and M. Raynaud, {\it N\'eron models}, Ergebnisse der Math. und ihrer Grenzgebiete (3), Vol. {\bf 21}, Springer-Verlag, Berlin, 1990.
\item{[2]} N. Bourbaki, {\it Lie groupes and Lie algebras}, Chapters {\bf 4--6}, Elements of Mathematics (Berlin), Springer-Verlag, Berlin, 2002.
\item{[3]} P. Deligne, {\it Travaux de Shimura}, S\'em. Bourbaki, Exp. no. 389, Lecture Notes in Math., Vol. {\bf 244},  123--165, Springer-Verlag, Berlin, 1971.
\item{[4]} P. Deligne, {\it Vari\'et\'es de Shimura: interpr\'etation modulaire, et techniques de construction de mod\`eles canoniques}, Automorphic forms, representations and $L$-functions (Oregon State Univ., Corvallis, OR, 1977), Part 2,   247--289, Proc. Sympos. Pure Math., Vol. {\bf 33}, Amer. Math. Soc., Providence, RI, 1979.
\item{[5]} P. Deligne, {\it Hodge cycles on abelian varieties}, Hodge cycles, motives, and Shimura varieties,
Lecture Notes in Math., Vol. {\bf 900},  9--100, Springer-Verlag, Berlin-New York, 1982.
\item{[6]} M. Demazure, A. Grothendieck, et al. {\it Sch\'emas en groupes, Vols. {\bf II--III}}, S\'eminaire de G\'eom\'etrie Alg\'ebrique du Bois Marie 1962/64 (SGA 3), Lecture Notes in Math., Vols. {\bf 152--153}, Springer-Verlag, Berlin-New York, 1970. 
\item{[7]} G. Harder, {\it \"Uber die Galoiskohomologie halbeinfacher Matrizengruppen II}, Math. Z. {\bf 92} (1966),  396--415.  
\item{[8]} S. Helgason, {\it Differential geometry, Lie groups, and symmetric spaces}, Pure and Applied Mathematics, Vol. {\bf 80}, Academic Press, Inc. [Harcourt Brace Jovanovich, Publishers], New York-London, 1978.
\item{[9]} R. E. Kottwitz, {\it Points on some Shimura varieties over finite fields}, J. Amer. Math. Soc. {\bf 5} (1992), no. 2,  373--444.
\item{[10]} J. S. Milne, {\it The points on a Shimura variety modulo a prime of good reduction}, The Zeta functions of Picard modular surfaces,  153--255, Univ. Montr\'eal, Montreal, QC, 1992.
\item{[11]} J. S. Milne, {\it Shimura varieties and motives}, Motives (Seattle, WA, 1991),  447--523, Proc. Sympos. Pure Math., Vol. {\bf 55}, Part 2, Amer. Math. Soc., Providence, RI, 1994.
\item{[12]} D. Mumford, J. Fogarty, and F. Kirwan, {\it Geometric invariant theory. Third edition}, Ergebnisse der Math. und ihrer Grenzgebiete (2), Vol. {\bf 34}, Springer-Verlag, Berlin, 1994.
\item{[13]} J.-P. Serre, {\it Groupes alg\'ebriques associ\'es aux modules de Hodge--Tate}, Journ\'ees de G\'eom. Alg. de Rennes (Rennes, 1978), Vol. III,  155--188, Ast\'erisque, Vol. {\bf 65}, Soc. Math. de France, Paris, 1979.
\item{[14]} G. Shimura, {\it On analytic families of polarized abelian varieties and automorphic functions}, Ann. of Math. (2) {\bf 78} (1963), no. 1,  149--192.\item{[15]} J. Tate, {\it $p$-divisible groups}, Proceedings of a conference on local fields (Driebergen, 1966),  158--183, Springer-Verlag, Berlin, 1967.
\item{[16]} J. Tits, {\it Reductive groups over local fields}, Automorphic forms, representations and $L$-functions (Oregon State Univ., Corvallis, OR, 1977), Part 1,   29--69, Proc. Sympos. Pure Math., Vol. {\bf 33}, Amer. Math. Soc., Providence, RI, 1979.
\item{[17]} A. Vasiu, {\it Integral canonical models of Shimura varieties of preabelian type}, Asian J. Math. {\bf 3} (1999), no. 2,  401--518.
\item{[18]} A. Vasiu, {\it A purity theorem for abelian schemes},
Michigan Math. J. {\bf 52} (2004), no. 1, 71--81.
\item{[19]} A. Vasiu, {\it Integral models in unramified mixed characteristic (0,2) of hermitian orthogonal Shimura varieties of PEL type, Part I}, to appear in J. Ramanujan Math. Soc., available at math.NT/0307205.
\item{[20]} A. Vasiu and T. Zink, {\it Purity results for $p$-divisible groups and abelian schemes over regular bases of mixed characteristic}, Doc. Math. {\bf 15} (2010), 571--599. 
\bigskip
\noindent
Adrian Vasiu, 

\noindent
\hbox{Department of Mathematical Sciences, Binghamton University,}

\noindent
\hbox{Binghamton, P. O. Box 6000, New York 13902-6000, U.S.A.}

\noindent
e-mail: adrian@math.binghamton.edu, fax: 1-607-777-2450, phone: 1-607-777-6036.

\end
We check that there exists a class $\gamma^\prime\in H^1(F,{\it Aut}^\prime(\Mg^{\rm ad}_F))$ that defines the form $\Md_{n,\tau,F}$ of $\Mg$ (i.e., whose image in $H^1(F,{\it Aut}(\Mg^{\rm ad}_F))$ is $\gamma$). For this we can assume that $n=4$. But if $n=4$, from [13, Cor. 2 of p. 182] we easily get that the image of $\gamma$ in $H^1(F,{\bf E}_F)$ is trivial and this implies that the class $\gamma^\prime$ exists. 
\smallskip